\documentclass[smallextended]{svjour3}
\usepackage[T1]{fontenc}
\usepackage{amsmath}
\usepackage{amsfonts}
\usepackage{amsxtra}
\usepackage{esint}
\usepackage[dvips]{epsfig}
\usepackage[utf8]{inputenc}
\usepackage{color}
\usepackage{bbm}
\usepackage{amssymb, latexsym}
\usepackage{booktabs}
\usepackage{setspace}
\usepackage{caption}
\usepackage{subfigure}
\usepackage{ulem}

\makeatother

\newtheorem{testproblem}{Test Problem}

\newcommand{\norm}[1]{\left\Vert#1\right\Vert}
\newcommand{\abs}[1]{\left\vert#1\right\vert}
\newcommand{\set}[1]{\left\{#1\right\}}
\newcommand{\real}{\mathbb R}
\newcommand{\setR}{\mathbb R}
\newcommand{\setN}{\mathbb N}
\newcommand{\DivG}{\nabla_\Gamma \cdot}
\newcommand{\diam}{\operatorname{diam}}
\newcommand{\landau}{\mathcal O}
\newcommand{\T}{\mathcal T _h}
\newcommand{\Tf}{\bar {\mathcal{ T}} _h}
\newcommand{\LF}{\text{LF}}

\title{Geometric Error of Finite Volume Schemes for Conservation Laws on Evolving Surfaces}
\author{Jan Giesselmann and Thomas M\"uller}
\titlerunning{Geometric Error of FV Schemes for Conservation Laws on Evolving Surfaces}
\institute{J. Giesselmann
 \at Weierstrass Institute,\\
     Mohrenstr. 39, D-10117 Berlin, Germany\\
     \email{giesselm@wias-berlin.de}
 \and 
T. M\"uller
 \at Abteilung f\"ur Angewandte Mathematik, Universit\"at Freiburg,\\
      Hermann-Herder-Str. 10, D-79104 Freiburg, Germany\\
      \email{mueller@mathematik.uni-freiburg.de}}

\begin{document}

\maketitle

\begin{abstract}

This paper studies finite volume schemes for scalar hyperbolic conservation laws on evolving hypersurfaces of $\mathbb{R}^3$. 
We compare theoretical schemes assuming knowledge of all geometric quantities to (practical) schemes defined on moving polyhedra approximating the surface.
For the former schemes error estimates have already been proven, but the implementation of such schemes is not feasible for complex geometries.
The latter schemes, in contrast, only require (easily) computable geometric quantities and are thus more useful for actual computations. 
We prove that the difference between approximate solutions defined by the respective families of schemes is of the order of the mesh width. 
In particular, the practical scheme converges to the entropy solution with the same rate as the theoretical one.
Numerical experiments show that the proven order of convergence is optimal. 

\end{abstract}

\section{Introduction}

Hyperbolic conservation laws serve as models for a wide variety of applications in continuum dynamics.
In many applications the physical domains of these problems
are stationary or moving hypersurfaces. Examples of the former are in particular geophysical problems \cite {WDHJS92}
and magnetohydrodynamics in the tachocline of the sun \cite{Gil00,SBG01}.
Examples of the latter include transport processes on cell surfaces \cite{RS05}, surfactant flow
on interfaces in multiphase flow \cite{BPS05} and petrol flow on a time dependent water surface.
There are several recent approaches to the numerical computation of such equations. 
Numerical schemes for the shallow water equations on a rotating sphere can be found in \cite{CHL06,Gir06,Ros06}.
For the simulation of surfactant flow on interfaces we refer to \cite{AB09,BS10,JL04}.
As we are interested in numerical analysis we focus on nonlinear scalar conservation laws as a model for these systems.
The intense study of conservation laws posed on fixed Riemannian manifolds started within the last years.
There are results on well-posedness \cite{BL06,DKM13,LM13} of the differential equations and on
the convergence of appropriate finite volume schemes \cite{ABL05,Gie09,GW12,LON09}.
For recent developments on finite volume schemes for parabolic equations we refer to \cite{LNR11}.

In the previous error analysis for finite volume schemes approximating nonlinear conservation laws on manifolds
 the schemes were defined on curved elements lying on the curved surface and
 it was assumed that geometric quantities like lengths, areas and conormals are known exactly.
While this is a reasonable assumption for schemes defined on general Riemannian manifolds or even more
 general structures \cite{ALN11,LO08.2} with no ambient space, most engineering applications involve equations
 on hypersurfaces of $\real^3$ and one aims at computing the geometry with the least effort.
This is in particular important for moving surfaces where the geometric quantities have to be computed in each time step.
Now the question arises to which extent an approximation of the geometry influences the order of convergence
 of the scheme.
 We will treat this question in the ``embedded`` case where an explicit embedding of the surface under consideration into Euclidean space is known.
It is an interesting question for future studies whether our analysis can be extended to errors arising from the discretisation of the geometry of the underlying space 
in an ''invariant'' description,
like the very general one analysed in \cite{ALN11,LO08.2}.

We consider the following initial value problem, posed on a family of closed, smooth hypersurfaces
 $\Gamma=\Gamma(t) \subset \real^3$. For a derivation cf. \cite{DE07.2,DKM13,Sto90}.
For some $T>0$, find $u: G_T := \bigcup_{t \in [0,T]}\Gamma(t) \times \{t\} \rightarrow\real$ with
     \begin{align}
        \dot{u} + u \DivG v+  \DivG  f(u,\cdot,\cdot)  &= 0  & &\text{in } G_T, \label{eq:ConLaw}\\
                                     u(\cdot,0)&=u_0 & &\text{on } \Gamma(0),\label{eq:IVConLaw}
     \end{align}
where 
$v$ is the
 velocity of the material points of the surface 
and $u_0:\Gamma(0)\rightarrow\real$ are initial data. 
For every $\bar u\in \real,\ t\in[0,T]$ the flux $f(\bar u,\cdot,t)$  is a smooth vector field tangential to $\Gamma(t),$
which depends Lipschitz on $\bar u$ and smoothly on $t$.
Moreover, we impose the following growth condition
\begin{equation}\label{eq:growth} |\nabla_{\Gamma} \cdot f (\bar u,x,t) | \leq c + c |\bar u| \quad \forall \, \bar u \in \setR , (x,t) \in G_T\end{equation} 
for some constant $c>0$.
By $\dot{u}$ we denote the material derivative of $u$ which is given by
\[ \dot{u}(\Phi_t(x),t) := \frac{d}{dt} u(\Phi_t(x),t),\]
where $\Phi_t : \Gamma(0) \rightarrow \Gamma(t)$
is a family of diffeomorphisms depending smoothly on $t$, such that  $\Phi_0$ is the identity on $\Gamma(0).$
Obviously this excludes changes of the topology of $\Gamma.$
We will assume that the movement of the surface and also the family $\Phi_t$ is prescribed.
A main result of this paper is a bound for the difference between two approximations of $u$. In particular,
we will give an estimate for the difference between the flat approximate and the curved approximate solution.
By curved approximate solution we refer to a numerical solution given by a finite volume scheme defined on the curved surface,
cf. Section \ref{subsec:fv_curved}, and
by flat approximate solution we refer to a numerical solution given by a finite volume scheme defined on a polyhedron approximating the surface,
cf. Section \ref{subsec:fv_flat}.
We will see that the arising geometry errors can be neglected compared to the error between the curved approximate solution 
 and the exact solution, i.e. both approximate solutions converge to the entropy solution with the same convergence rate.
We will present numerical examples showing that the proven convergence rate is optimal under the assumptions for the numerical analysis.
However, for most numerical experiments we observe higher orders of convergence.

Our analysis also indicates that the geometry error poses an obstacle to the construction 
of higher order schemes. To this end we perform numerical experiments underlining in which manner
the order of convergence of the higher order scheme is restricted by the approximation of the geometry.
This shows that to obtain higher order convergence
also the geometry of the manifold has to be approximated more accurately, cf. \cite{Dem09} in a finite element context.

The outline of this paper is as follows. In Section \ref{sec:fvs} we review the definition of finite volume schemes on 
moving curved surfaces and define finite volume schemes on moving polyhedra approximating the surfaces.
The approximation errors for geometric quantities are established in Section \ref{sec:geom}.
Section \ref{sec:mr} is devoted to estimating the difference between the curved and the flat approximate solution.
Finally, numerical experiments are given in Section \ref{sec:numerics}.

\section{The Finite Volume Schemes}\label{sec:fvs}
This section is devoted to the construction of a family of triangulations $\T(t)$ of the surfaces suitably linked to 
 polyhedral approximations $\Gamma_h(t)$ of the surfaces. Afterwards we will recall the definition
 of a finite volume scheme on $\T(t)$ which was considered in the hitherto error analysis 
and define a finite volume scheme on $\Gamma_h(t)$ which is an algorithm only relying on
 easily computable quantities. 
 We would like to point out that our finite volume scheme is applicable to closed, smooth hypersurfaces of arbitrary geometry,
 which additionally may evolve in time. In the case of simpler geometries of special interest,
 e.g. considering conservation laws on $\mathbb{S}^2$, one can make use of the additional structures.
 In \cite{BFL09} for instance, studying scalar conservation laws on $\mathbb S^2$, the special structure
 of this setting is exploited by considering a  longitude-latitude grid
 which allows the exact computation of cell areas and edge lengths and the design of a Godunov-type finite volume scheme
 based on the dimension-wise solution of one-dimensional Riemann problems for a class of analytic flux functions.
 Additionally, the finite volume scheme in \cite{BFL09}
is geometry-compatible in the sense that for a divergence-free flux the numerical fluxes are (discretely) divergence-free, as well.
Other approaches employing the special knowledge available for the sphere include logically rectangular grids developed in \cite{CHL06}
and grids in which all edges are geodesic arcs in \cite{Gir06}.

 We mention that our triangulation as well as the definition of the finite volume
scheme on $\Gamma_h$ is in the same spirit as the one from \cite{LNR11} which was developped for
the diffusion equation on evolving surfaces.

\subsection{Triangulation}
We start by mentioning that there are neighbourhoods $\mathcal N(t)\subset \setR^3$ of $\Gamma(t)$
 such that for every $x\in \mathcal N(t)$ there is a unique
point $a(x,t)\in \Gamma(t)$ such that
\begin{equation}
x = a(x,t) + d(x,t) \nu_{\Gamma(t)}(a(x,t)),
 \label{eq:projection}
\end{equation}
where $d(\cdot,t)$ denotes the signed distance function to $\Gamma(t)$ and $\nu_{\Gamma(t)}(a(x,t))$ 
the unit normal vector to $\Gamma(t)$  pointing towards the non-compact component of $\setR^3\setminus \Gamma(t)$.
 See \cite{DE07} for example.

Let us choose a polyhedral surface $\Gamma_h(0) \subset \mathcal N(0)$ which consists of flat triangles
 such that the vertices 
of $\Gamma_h(0)$ lie on $\Gamma(0),$ and $h$ is the length of the longest edge of $\Gamma_h(0).$
 In addition we impose that the restriction of  $a(\cdot,0)|_{\Gamma_h(0)}: \Gamma_h(0) \rightarrow \Gamma(0)$ is one-to-one.
We define $\Gamma_h(t)$ as the polyhedral surface that is constructed by moving the vertices of $\Gamma_h(0)$ via the 
diffeomorphism $\Phi_t$ and connecting them  with straight lines such that all triangulations share 
the same grid topology.
A triangulation $\Tf(t)$ of $\Gamma_h(t)$ is automatically given by the decomposition into closed faces. 
We define the triangulation $\T(t)$ on $\Gamma(t)$ as the image of $\Tf(t)$ under $a(\cdot,t)|_{\Gamma_h(t)}.$
 We will denote the closed
curved cells with $K(t)$ and the closed curved faces with $e(t)$. 
A flat quantity corresponding to some curved quantity is denoted by the same letter and a bar, 
e.g. let $e(t) \subset \Gamma(t)$ be a curved face then $\bar e(t) = (a(\cdot,t)|_{\Gamma_h(t)})^{-1}(e(t)).$
In order to reflect the fact that all triangulations share the same grid topology we introduce the following 
notation. We denote by $K$ the family of all (closed) curved triangles relating to the same (closed) triangle
 $\bar K(0)$ on $\Gamma_h(0).$ We do the same for $e, \bar K, \bar e.$ Analogously by $\T$ we denote 
the family of such families
of triangles $K.$

For later use we state the following Lemma summarizing geometric properties, whose derivation can be found in \cite{DE07}.
\begin{lemma}
Let $\Gamma_h(t)$ be a polyhedral approximation of $\Gamma(t)$ as described above then there exists $C=C(T)$ such that 
for all $t \in [0,T]$ 
\begin{enumerate}
\item $\nu_{\Gamma(t)} = \nabla d(\cdot, t),$
\item $\| d(\cdot,t)|_{\Gamma_h(t)} \|_{L^\infty(\Gamma_h(t))} \leq C h^2$ .
\end{enumerate}
\label{lemma:dziuk-elliott}
\end{lemma}

We will use the following notation. By $h_{K(t)} := \diam(K(t))$ we denote the diameter of each cell, furthermore
 $h:= \max_{t \in [0,T]} \max_{K(t)}h_{K(t)}$ and $\abs{K(t)},$ $\abs{\partial K(t)}$ are the Hausdorff measures 
of $K(t)$ and the boundary of $K(t)$ respectively. 
When we write $e(t)\subset \partial K(t)$ we mean $e(t)$ to be a face of $ K(t)$.

We need to impose the following assumption uniformly on all triangulations $\Tf(t).$
 There is a constant number $\alpha >0$ such that for each flat cell $\bar K(t) \in  \Tf(t)$
 we have
\begin{align}
\begin{split}
 \alpha h^2 \leq \abs{\bar K(t)},  \\
 \alpha \abs{\partial \bar K(t)} \leq h.
\end{split}
\label{ineq:angle-assumption}
\end{align}
Later on, we will see that \eqref{ineq:angle-assumption} implies the respective estimate 
for the curved triangulation, cf. Remark \ref{rem:angles}.
A consequence of \eqref{ineq:angle-assumption} is that $2 \alpha^2 h$ is a lower bound of the 
radius of the inner circle of $\bar K(t)$, which implies that the sizes of the angles in $\bar K(t)$ are
 bounded from below.
Furthermore we denote by $\kappa(x,t)$  the supremum of the spectral
 norm of $\nabla \nu_{\Gamma(t)}(x),$ i.e.  $\kappa$ is a bound on the eigencurvatures.
By straightforward continuity and compactness arguments $\kappa$ is uniformly bounded in space and time.

\subsection{The Finite Volume Scheme on Curved Elements}\label{subsec:fv_curved}
In this section we will briefly review the notion of finite volume schemes on moving  curved surfaces.
We consider a sequence of times $0=t_0< t_1 < t_2 < \dots$ and set $I_n:= [t_n,t_{n+1}].$ Moreover we assign
 to each $n \in \setN$ and $K \in \T$
 the term $u^n_K$ approximating the mean value of $u$ on $\bigcup_{t \in  I_n} K(t)\times \{t\}$ and to
 each $K \in \T$ and face $e \subset \partial K$ a numerical flux function $f^n_{K,e} : \real^2 \rightarrow \real$,
 which should approximate 
\begin{equation}\label{curv_num_flux_approximates}
\fint_{I_n}
\fint_{e(t)} \langle f(u(x,t),x,t), \mu_{K(t),e(t)}(x) \rangle \, de(t) \, dt,
\end{equation}
where $de(t)$ is the line element, $\mu_{K(t),e(t)}(x) $ is the unit conormal to $e(t)$ pointing outwards from $K(t)$ and
$\langle \cdot, \cdot \rangle $ is the standard Euclidean inner product.
 Please note that $\mu_{K(t),e(t)}(x) $ is tangential to $\Gamma(t).$
Then the finite volume scheme is given by
\begin{align}\begin{split}\label{eq:FVcurv}
u^0_K &:= \fint_{K(0)} u_0(x) d\Gamma(0),\\
u^{n+1}_K&:= \frac{|K(t_n)|}{|K(t_{n+1})|} u^n_K - \frac{\abs{I_n}}{|K(t_{n+1})|}
\sum_{e \subset \partial K} |e(t_n)| f^n_{K,e}(u_K^n,u_{K_e}^n),\\
u^h(x,t)&:= u^n_K \quad \text{ for }  t \in [t_n,t_{n+1}),\ x \in K(t),\end{split}
\end{align}
where $K_e$ denotes the cell sharing 
 face $e$ with $K$ and $d\Gamma(0)$ is the surface element.
As usual in corresponding convergence analysis \cite{Gie09,LON09} we assume that the used numerical fluxes are consistent, i.e.
\begin{equation}\label{def:consistency}
|e(t_n)| f^n_{K,e} (u,u)= \fint_{I_n} \int_{e(t)} \langle f(u,x , t),\mu_{K(t),e(t)}(x) \rangle de(t) dt \quad \forall u \in \setR,
\end{equation}
 conservative, i.e.
\begin{equation}\label{def:conserv}
 f^n_{K,e} (u,v)= -f^n_{K_e,e} (v,u) \quad \forall u,v \in \setR,
\end{equation}
 monotone, i.e.
\begin{equation}\label{def:mono}
 \frac{d}{du} f^n_{K,e} (u,v) \geq 0, \quad \frac{d}{dv} f^n_{K,e} (u,v) \leq 0\quad \forall u,v \in \setR,
\end{equation}
and uniformly Lipschitz continuous.
Let $L$  denote the Lipschitz constant of the numerical fluxes, then additionally the CFL condition
\begin{equation}\label{cfl}t_{n+1}-t_n \leq \frac{\alpha^2 h}{8L}\end{equation}
has to be imposed to ensure stability, i.e. 
Lemmas \ref{lem:curvedstability} and \ref{lem:contraction}.
As an example for numerical flux functions satisfying these conditions we introduce Lax-Friedrichs fluxes
\begin{multline}
\label{curv_num_flux}  ^{\LF}f^n_{ K,e} (u,v) :=  \fint_{I_n} 
 \frac{1}{2|e(t_n)|}\int_{e(t)}
 \langle  f(u,x,t)+ f(v,x,t), \mu_{ K(t), e(t)}(x) \rangle\, de(t)\, dt\\ + \lambda_n (u-v),
\end{multline}
where $\lambda_n\geq \frac{1}{2}\|\partial_u f \|_\infty \frac{\max_{t \in I_n} |e(t)|}{|e(t_n)|}$ is an artificial viscosity coefficient ensuring the monotonicity of 
$f^n_{K,e}$ and stabilizing the scheme.

\subsection{The Finite Volume Scheme on Flat Elements}
\label{subsec:fv_flat}
In this section we define finite volume schemes on $\Tf$ which are
in the same spirit as \eqref{eq:FVcurv} but 
only rely on easily accessible geometrical information.
We want to point out that the calculation of areas and lengths is straightforward for flat elements. 
As well, the approximation of integrals can be achieved using quadrature formulas by mapping cells and edges
 to a standard triangle and  the unit interval, respectively, using affine linear maps.
In this fashion we obtain for every time $t \in [0,T]$ quadrature operators $Q_{\bar K(t)}: C^0(\bar K(t)) \rightarrow \setR,$ and 
$ Q_{\bar e(t)}: C^0(\bar e(t)) \rightarrow \setR$
 of order $p_1,p_2 \geq 1,$ respectively. In addition 
 for any compact interval $I \subset [0,T]$  the term
$Q_I: C^0(\bar I) \rightarrow \setR$
denotes a quadrature operator of order $p_3\geq 1.$
For Lipschitz continuous numerical flux functions $\bar f_{\bar K,\bar e}^n:\setR\times\setR\rightarrow\setR$ we define the
finite volume scheme on flat elements according to
\begin{align}\begin{split}
\label{eq:FVflat} 
\bar u^0_{\bar K} &:= \frac{1}{|\bar K(0)|} Q_{\bar K(0)} (u_0(a(\cdot,0))),\\
 \bar u^{n+1}_{\bar K}&:=\frac{|\bar K(t_n)|}{|\bar K(t_{n+1})|} \bar u^n_{\bar K} 
- \frac{|I_n|}{|\bar K(t_{n+1})|}\sum_{\bar e \subset\partial \bar K} |\bar e(t_n)| \bar f^n_{\bar K,\bar e}(\bar u_{\bar K}^n,\bar u_{\bar K_{\bar e}}^n),\\
\bar u^h(x,t)&:= \bar u^n_{\bar K}, \quad \text{ for }  t \in [t_n,t_{n+1}),\ x \in  K(t).
\end{split}\end{align}
Note that by \eqref{eq:FVflat}$_3$ the function $\bar u^h$ is defined on $G_T$.
For the numerical analysis we need to impose the following  estimate for the (geometric) error between the numerical fluxes $f_{K,e}^n$ and $\bar f_{\bar K,\bar e}^n$:
\begin{align}\label{eq:fluxerror}
\left| f^n_{K,e} (u,v) - \bar f^n_{\bar K, \bar e} (u,v)\right|\leq Ch^2 
\qquad \forall \ (u,v) \in \mathcal{K} ,\ K\in\mathcal T_h,\ e\subset \partial K,
\end{align}
where $\mathcal{K}$ is a compact subset of $\setR^2$ and $C$ a constant depending only on $G_T$ and $\mathcal{K}$.

As an example for easily computable numerical flux functions for the flat scheme we define
Lax-Friedrichs flux functions below.
We will see in Lemma \ref{lem:flux_est} that assumption \eqref{eq:fluxerror} is valid for them.
Before we can use the quadrature operators to define the numerical fluxes we need to determine the "discrete" conormals.
To each flat triangle $\bar K(t)$ we fix a unit normal $\bar \nu_{\bar K(t)}$ by imposing 
\begin{equation}\label{**} 
\langle\bar \nu_{\bar K(t)}, \nu_{\Gamma(t)}(y)\rangle >0,
\end{equation}
where $y$ is the barycentre of $K(t).$
We will see in Lemma \ref{lem:angles} that $\bar \nu_{\bar K(t)}$ converges 
to $\nu_{\Gamma(t)}(y)$ for $h \rightarrow 0.$
 To each face $\bar e(t)$ and adjacent cell $\bar K(t)$ 
 there is a unique unit tangent vector $\bar {\bf t}_{\bar K(t),\bar e(t)}$ such that
$\bar \nu_{\bar K(t)} \times \bar {\bf t}_{\bar K(t), \bar e(t)} $ is a conormal to $\bar e(t)$ pointing outward
 from $\bar K(t).$
Hence  this vector product is one candidate for 
$\bar \mu_{\bar K(t),\bar e(t)}.$
  However in general 
\begin{equation}\label{lack_of_conserv}
\bar \nu_{\bar K(t)} \times \bar {\bf t}_{\bar K(t),\bar e(t)} \not= \pm (\bar \nu_{{\bar K}_{\bar e}(t)} \times 
\bar {\bf t}_{\bar K_{\bar e}(t),\bar e(t)})
\end{equation}
such that  a choice like
\[ \bar \mu_{\bar K(t),\bar e(t)} =\bar \nu_{\bar K(t)} \times \bar {\bf t}_{\bar K(t),\bar e(t)}\]
 would lead to a loss of conservativity of the resulting numerical fluxes.
Therefore we choose
\[ \bar \mu_{\bar K(t), \bar e(t)}:= \frac{1}{2} \left( \bar \nu_{\bar K(t)} \times \bar {\bf t}_{\bar K(t),\bar e(t)}
 + \bar \nu_{{\bar K}_{\bar e}(t)} \times \bar {\bf t}_{\bar K(t),\bar e(t)}\right).\]
We define numerical Lax-Friedrichs fluxes by
\begin{align}\begin{split}
\label{eq:LFflat}
^\LF\bar f_{\bar K,\bar e}^n (u,v) &:= \frac{1}{|I_n|} Q_{I_n}\left[\frac{1}{2\abs{\bar e(t_n)}} 
Q_{\bar e(\cdot)}\left( \langle f(u,\cdot,\cdot)+ f(v,\cdot,\cdot),
  \bar \mu_{\bar K(\cdot), \bar e(\cdot)}\rangle \right)\right] \\
&\qquad \qquad \qquad \qquad  \qquad \qquad \qquad \qquad \qquad  \qquad \quad + \lambda (u-v)
\end{split}
\end{align}
for some sufficiently large $\lambda \geq 0$
and $f$ being smoothly extended from $G_T$ to the whole of $\bigcup_{t\in [0,T]} \mathcal{N}(t) \times \{t\}$.
Note, that here and in the following the quadrature operator $Q_{I_n}$ is applied to the time dependence
while $Q_{\bar e(\cdot)}$ is applied to the space dependence.
In particular, for each time quadrature point $\tau_i$, specified by $Q_{I_n}$, the 
space quadrature points lie on the (moving) flat edge $\bar e (\tau_i)$.

\section{Geometrical Estimates}\label{sec:geom}
In this section we derive estimates for the approximation errors of the geometric quantities.
Throughout this section we suppress the time dependence of all quantities.
 All the estimates can be derived uniformly in time.
To obtain the geometrical estimates, we introduce the following lift operator.

\begin{definition}
Let $ \bar U \subset \Gamma_h$ and $\bar g$ a function on $\bar U$ then we define a function $\bar g^l$ on
 $a|_{\Gamma_h}(\bar U)$ as 
\[ \bar g^l= \bar g \circ a|_{\Gamma_h}^{-1}.\]
Similarly we define the inverse of this lift operator by
\[ g^{-l}= g \circ a|_{\Gamma_h}\]
for a function $g$ defined on some $ U \subset \Gamma$.
\end{definition}

We begin our investigation with the differences between the normal
vectors of the flat and curved elements.

\begin{lemma}\label{lem:angles}
There is a constant $C$ such that for all flat cells $\bar K$ 
and every 
$y \in \bar K$ we have
\begin{align}
\label{norm_est} \left\| \nu_\Gamma^{-l}(y) - \bar \nu_{\bar K} \right\| &\leq Ch.
\end{align}
The constant $C$ depends on derivatives of $d,$ in particular on $\kappa.$
\end{lemma}

\proof
Without loss of generality we can assume that
$\bar K $ is a subset of $\{ (x,y,0) \in \setR^3 \, | \, y<0\}$ 
such that $\bar e= \{ (s,0,0) \in \setR^3 \, |\, s \in [0,h_{\bar e}]\}$ is one of its faces
and $(\nu_\Gamma)^{-l}_3(y)>0$ for some $y \in \bar K$.
We start by showing that there exists some constant $C>0$ such that
\begin{equation}\label{nuk_est} |(\nu_\Gamma)_i|\leq Ch , \quad \text{for } i=1,2.\end{equation}
We recall that $\nu_\Gamma=\nabla d,$ where $d$ is the signed distance function to $\Gamma.$
 As the vertices of $\Gamma_h$ lie on $\Gamma$ we know that, if we denote the third vertex by $(x,y,0) \in \bar K,$ it holds 
\[ d(0,0,0)=0,\ d(h_{\bar e},0,0)=0, \ d(x,y,0)=0.\]
Hence, the directional derivatives of $d$ with respect to $(x,y,0)$ and $(1,0,0)$
 need to vanish somewhere in $\bar K.$ 
Thus their absolute value is of order $\landau(h)$ on $\bar K.$
 Due to the angle condition \eqref{ineq:angle-assumption}
 an analogous inequality also holds for the directional derivative of $d$ with respect to
 $(0,1,0)$. As the directional derivative of $d$ with respect to $(1,0,0),\ (0,1,0)$ coincides with 
$(\nu_\Gamma)_1,\ (\nu_\Gamma)_2,$ respectively, this proves \eqref{nuk_est}.
This immediately implies $(\nu_\Gamma)_3 =\pm \sqrt{1- \landau(h^2)} = \pm 1 + \landau(h^2).$ 
By assumption  $(\nu_\Gamma)_3 = 1 + \landau(h^2)$ everywhere and by \eqref{**} we have
 $\bar \nu_{\bar K}=(0,0,1) $ which proves \eqref{norm_est}.
\qed

\begin{lemma}
 \label{lem:eK}
For the difference between the length of a curved edge $e$ and the corresponding flat edge $\bar e$ we have
\begin{equation} 
 \abs{\frac{\abs{e}}{\abs{\bar e }} - 1  } \leq C h^2,
\label{eq:error-e}
\end{equation}
and for the difference between the area of a curved cell $K$ and the corresponding
flat cell $\bar K$  we have
\begin{equation}
 \abs{\frac{\abs{K}}{\abs{\bar K }} - 1  } \leq C h^2,
\label{eq:error-K}
\end{equation}
where $C$ does not depend on $h$ but on $\kappa.$

Furthermore let $c_e$ be the parametrization of $e$ over $\bar e$ given by $a|_{\bar e}$ then we have
\begin{equation}\label{le:normc}
\abs{ \norm{c_e^\prime(s)}- 1}\leq C h^2.
\end{equation}
\end{lemma}

\proof
We assume without loss of generality that $\bar K \subset \{ (x,y,0) \in \setR^3 \, | \, y<0\}$. For small
 enough $h$ we can parametrize the curved cell $K$
according to \eqref{eq:projection} by a parametrization $c = a|_{\bar K} : \bar K \rightarrow K \subset \real^3$ with
\begin{equation*}
c(x_1,x_2) = (x_1,x_2,0) - d(x_1,x_2,0) \nu_\Gamma(c(x_1,x_2)),
 \label{eq:c-K}
\end{equation*}
where we suppressed the third coordinate in $\bar K$.
The ratio of volume elements of $K$ and $\bar K$ with respect to the parametrization $c$ is given by
\begin{equation*}
 \sqrt{\abs {g}}:=\sqrt{\det(g)},
\end{equation*}
where the matrix $g$ is defined by
\begin{equation*}
 g = \left(g_{ij}\right)_{1\leq i,j \leq 2} : =
 \left(\left\langle \partial_i c, \partial_j c\right\rangle\right)_{1\leq i,j \leq 2}.
\end{equation*}
For the parametrization $c$ of $K$ we have
\begin{equation*}
 \partial_i c = e_i - \left\langle \nabla d , e_i \right\rangle \nu_\Gamma\circ c - d \;
 \partial_i c\; \left( \nabla \nu_\Gamma\right)^T \circ c  \quad\text{for } i=1,2,
\end{equation*}
where $e_i$ denotes the $i$-th standard unit vector.
Due to the bounded curvature of $\Gamma$ and Lemma \ref{lemma:dziuk-elliott} we can show that
\begin{equation}\label{eq:partialc}
 \partial_i c = e_i - ((\nu_\Gamma)_i \nu_\Gamma)\circ c +  \mathcal  O ( h^2) \quad\text{for } i=1,2.
\end{equation}
Applying \eqref{norm_est}  we see that 
\begin{equation*}
 \nu_\Gamma= \pm(0,0,1) +  \mathcal O( h) \text{ and } 
\left\langle e_i,\nu_\Gamma\right\rangle =(\nu_\Gamma)_i=  \mathcal O (h)\quad\text{for } i=1,2.
\end{equation*}
Thus, for the matrix $g$ we have
\begin{equation*}
 g = \begin{pmatrix}
      1+ \mathcal O ( h^2) &  \mathcal O ( h^2)\\
       \mathcal O ( h^2) & 1+ \mathcal O ( h^2)
     \end{pmatrix}
\end{equation*}
which implies for the volume element 
\begin{equation}
 dK= \sqrt{\abs {g}}d\bar K = \sqrt{1 +  \mathcal O (h^2)}d\bar K = d \bar K +  \mathcal O (h^2)d\bar K.
\end{equation}
Therefore, we arrive at
\begin{align*}
 \abs{\abs{K} - \abs{\bar K } } = \abs{\int_{\bar K} \sqrt{\abs {g}}  - 1d\bar K }
\leq C \abs{\bar K } h^2
\end{align*}
for the error of the cell area which proves \eqref{eq:error-K}.

To prove \eqref{eq:error-e} and \eqref{le:normc} we consider without loss of generality an edge 
$\bar e = \set{(s,0,0)| 0\leq s \leq h_{\bar e}} \subset \partial\bar K$, where $h_{\bar e}$ denotes the length of $\bar e $.
Considering the derivation of the parametrization 
\begin{equation}
c_{ e}(s) = c(s,0) =  (s,0,0) - d(s,0,0) \nu_\Gamma(c_{ e}(s))
 \label{eq:c-e}
\end{equation}
of the curved edge $e$ and applying the same arguments as we used to prove \eqref{eq:error-K} completes the proof.
\qed

\begin{remark}
Let us note that an analogous estimate to \eqref{ineq:angle-assumption} for curved elements is an easy consequence
 of \eqref{ineq:angle-assumption}, \eqref{eq:error-e}, \eqref{eq:error-K} and the fact $|h_{\bar K}-h_K | \leq C h^2,$
 which is a consequence of Lemma \ref{lem:eK}.
\label{rem:angles}
\end{remark}

\begin{lemma}\label{lem:angles1}
There is a constant $C$ (depending on $\kappa$) such that for all flat cells $\bar K,$ all flat edges $\bar e  \subset  \partial \bar K$ and every $x \in \bar e$ we have
\begin{align}
\label{conorm_est1} \left|\langle \bar \mu_{\bar K,\bar e} , {\bf t}^{-l}(x) \rangle \right| &\leq C h^2,\\
\label{conorm_est2} \left|\langle \bar \mu_{\bar K,\bar e} , \nu_\Gamma^{-l}(x) \rangle \right| &\leq C h,\\
\label{conorm_est3} \left|\langle \bar \mu_{\bar K,\bar e} , \mu_{K,e}^{-l}(x) \rangle - 1\right| &\leq C h^2,
\end{align}
where ${\bf t}$ denotes a unit tangent vector to $e$. We want to point out that this estimate is independent of the sign 
of ${\bf t}.$
\end{lemma}
\proof
It is sufficient to show versions of \eqref{conorm_est1} - \eqref{conorm_est3}
 where $\bar \mu_{\bar K,\bar e}$ is substituted by
$ \bar \nu_{\bar K} \times \bar {\bf t}_{\bar K,\bar e}.$ 
Then analogous results for  $ \bar \nu_{\bar K_{\bar e}} \times \bar {\bf t}_{\bar K,\bar e}$ are immediate 
and therefore estimates \eqref{conorm_est1} - \eqref{conorm_est3} follow because $\bar \mu_{\bar K,\bar e}$ is the mean of the vectors
$ \bar \nu_{\bar K_{\bar e}} \times \bar {\bf t}_{\bar K,\bar e}$ and 
$ \bar \nu_{\bar K} \times \bar {\bf t}_{\bar K,\bar e}$.
Firstly, we address the proof of \eqref{conorm_est1}.
Let the same assumptions as in the proof of Lemma \ref{lem:angles} hold and in addition let
$\bar e$ be given by $ \{ (x,0,0) \in \setR^3 \, |\, x \in [0,h_{\bar e}]\}$. 
We obviously have \begin{equation}\label{bar_mu}  \bar \nu_{\bar K} \times \bar {\bf t}_{\bar K,\bar e}= (0,1,0).\end{equation}
Note that the assumptions of the proof of Lemma \ref{lem:eK} are satisfied. Hence we can use \eqref{eq:c-e},
i.e. the parametrization of $e$ 
given by $c_{e}$ satisfies
\begin{equation}
c_{e}'(s)= \left(  1 , 0 , 0 \right) - \nu_\Gamma(c_{e}(s)) (\nu_\Gamma)_1(c_{e}(s))+ \landau(h^2) ,
\end{equation}
and, by definition of  ${\bf t}$,
it holds ${\bf t}(c_e(s))= c_e'(s)/\|c_e'(s)\|.$
Hence, in view of \eqref{le:normc}
 we obtain
\begin{equation}\label{t}
{\bf t}^{-l}(x) = ( 1 ,0,0) - \nu_\Gamma(c_{e}(s)) (\nu_\Gamma)_1(c_{e}(s)) + \landau(h^2)
\end{equation}
for some $s\in [0,h_{\bar e}]$.
Combining \eqref{bar_mu} and \eqref{t} we find using \eqref{nuk_est}
\[\left|\langle \bar \nu_{\bar K} \times \bar {\bf t}_{\bar K,\bar e}, {\bf t}^{-l}(x) \rangle \right| 
= \left| (\nu_\Gamma)_2(c_{e}(s)) (\nu_\Gamma)_1(c_{e}(s)) \right| + \landau(h^2) \leq C h^2,\]
 which is \eqref{conorm_est1}. Concerning \eqref{conorm_est2},
\[ \left|\langle \bar \nu_{\bar K} \times \bar {\bf t}_{\bar K,\bar e} , \nu_\Gamma^{-l}(x) \rangle \right| \leq \left| (\nu_\Gamma^{-l})_2(x) \right| \leq Ch\]
holds because of  \eqref{bar_mu} and \eqref{nuk_est}.
Thus, it remains to show \eqref{conorm_est3}.
By definition ${\bf t}^{-l}(x), \nu_\Gamma^{-l}(x), \mu_{K,e}^{-l}(x)$ form an orthonormal 
basis of $\setR^3$ and the vector 
$\bar \nu_{\bar K} \times \bar {\bf t}_{\bar K,\bar e}$ 
is of unit length. This means that for every $\bar x$ in $\bar e$ there exist 
$b_1(\bar x),b_2(\bar x),b_3(\bar x) \in \setR$ satisfying $b_1^2(\bar x) + b_2^2(\bar x) + b_3^2(\bar x)=1$ such that
\begin{equation}\label{eq:bs} \bar \nu_{\bar K} \times \bar {\bf t}_{\bar K,\bar e} = b_1(\bar x) {\bf t}^{-l}(\bar x) 
 + b_2(\bar x) \nu_\Gamma^{-l}(\bar x) + b_3(\bar x) \mu_{K,e}^{-l}(\bar x).
\end{equation}

We know from \eqref{conorm_est1} and \eqref{conorm_est2} that $|b_1(\bar x)|,|b_2(\bar x)| \leq C h$ for some $C>0,$ 
which implies using Taylor expansion
\begin{equation}\label{eq:b3} b_3(\bar x) = \pm \sqrt{ 1 + \landau(h^2)} = \pm 1 + \landau(h^2).\end{equation}
Note that it only remains to show that in \eqref{eq:b3} the ``+'' holds. As $b_3$ depends continuously on $\bar x$ it is 
sufficient to find one $(\bar x_1,0,0)\in\bar K$ such that
$ b_3(\bar x_1) =  1 + \landau(h^2).$
To that end we consider the curve $\gamma(s):= c(\bar x_1, s)$ for $s<0$ and small $|s|$ where $c$ is the parametrization
of $K$ from Lemma \ref{lem:eK}.
As $\gamma$ is leaving $K$ through $e$ we have
\begin{equation}\label{eq:angle_sign}
 0 < \langle  \gamma'(0) , \mu_{K,e}( \gamma(0)) \rangle .
\end{equation}
Due to \eqref{eq:bs}, \eqref{eq:b3} and the fact that $\mu_{K,e}$ is of unit length we already know that
\begin{equation}\label{eq:muke}
 \mu_{K,e} \equiv \pm ( 0 , 1, 0 ) + \landau(h).
\end{equation}
Inserting \eqref{eq:partialc} for $i=2$ and \eqref{eq:muke} in \eqref{eq:angle_sign} leads, together with 
Lemma \ref{lem:angles}, to the ``+'' in \eqref{eq:b3}, which completes the proof.
\qed\\

\section{Estimating the Difference Between Both Schemes }
\label{sec:mr}

This section is devoted to establishing a bound for the difference between the curved and flat approximate solutions.
To start with we show that the Lax-Friedrichs numerical fluxes from Section \ref{sec:fvs} satisfy
assumption \eqref{eq:fluxerror}.

\begin{lemma}\label{lem:flux_est}
Let $\mathcal{K}$ be some compact subset of $\setR^2$.
Then there is a constant $C$ depending only on $G_T$ and $\mathcal{K}$ such that for the Lax-Friedrichs fluxes 
\eqref{curv_num_flux} and \eqref{eq:LFflat} with the same diffusion rate $\lambda$ the following inequality holds
\[\left| {^\LF}f^n_{K,e} (u,v) - {{^\LF}\bar f^n_{\bar K, \bar e}} (u,v)\right|\leq Ch^2 
\qquad \forall \ (u,v) \in \mathcal{K} ,\ K\in\mathcal T_h,\ e\subset \partial K.\]
\end{lemma}

\proof
We start by observing that the diffusive terms drop out, such that
\begin{align}
2 \left| {^\LF}f^n_{K,e} (u,v) - {^\LF}\bar f^n_{\bar K, \bar e} (u,v)\right| = 
|E^n_{K,e} (u) + E^n_{K,e} (v)|
\label{eq:j1}
\end{align}
with
\begin{align*}
\begin{split}
 E^n_{K,e} (u):=&\fint_{I_n} \frac{1}{|e(t_n)|} \int_{e(t)} \langle f(u,x,t), \mu_{K(t),e(t)}(x)\rangle de(t)\, dt\\
 &-\frac{1}{|I_n|} Q_{I_n}\left[ \frac{1}{|\bar e(t_n)|} Q_{\bar e(\cdot)}
 [\langle f(u,\cdot,\cdot), \bar \mu_{\bar K(\cdot),\bar e(\cdot)}\rangle ]\right].
\end{split}
\end{align*}
As $u$ and $v$ appear symmetrically
in \eqref{eq:j1}, we focus on the error analysis of only $E^n_{K,e} (u)$.

Addition of several zeros leads to
\begin{align}\label{eq:flux_est}
 |E^n_{K,e} (u)| = \left|\fint_{I_n}\frac{1}{|\bar e(t_n)|} \big(T_1 + T_2 + T_3 +T_4 +T_5\big) \; dt\right|
\end{align}
with
\begin{align*}\begin{split}
&T_1(t) := \frac{|\bar e(t_n)|}{|e(t_n)|}\int_{e(t)} \langle f(u,x,t), \mu_{K(t),e(t)}(x)\rangle de(t)  -
 \int_{\bar e(t)}\langle  f^{-l}(u,x,t), \mu_{K(t),e(t)}^{-l}(x)\rangle d \bar e(t), \\
&T_2(t):= \int_{\bar e(t)}\langle f^{-l}(u,x,t) ,\mu_{K(t),e(t)}^{-l}(x)\rangle d\bar e (t)
- \int_{\bar e(t)} \langle f^{-l}(u,x,t), \bar \mu_{\bar K(t),\bar e(t)}\rangle d\bar e(t),\\
&T_3(t):= \int_{\bar e(t)} 
\langle f^{-l}(u,x,t), \bar \mu_{\bar K(t),\bar e(t)}\rangle d\bar e(t) 
 - \int_{\bar e(t)}\langle f(u,x,t), \bar \mu_{\bar K(t),\bar e(t)}\rangle d \bar e(t), \\
&T_4(t):= \int_{\bar e(t)} 
\langle f(u,x,t) ,\bar \mu_{\bar K(t),\bar e(t)} \rangle d\bar e(t)
-\frac{1}{|I_n|}Q_{I_n}\left[ \int_{\bar e(\cdot)} 
\langle f(u,x,\cdot) ,\bar \mu_{\bar K(\cdot),\bar e(\cdot)} \rangle d\bar e(\cdot)  \right],\\
&T_5:=\frac{1}{|I_n|} Q_{I_n}\!\left[ \int_{\bar e(\cdot)} 
\langle f(u,x,\cdot) ,\bar \mu_{\bar K(\cdot),\bar e(\cdot)} \rangle d\bar e(\cdot)  
 -Q_{\bar e(\cdot)}\left[ \left\langle f(u,\cdot ,\cdot), \bar \mu_{\bar K(\cdot),\bar e(\cdot)}
\right\rangle\right] \right]\!.
\end{split}
\end{align*}
In the following we will
estimate the summands one by one.
First, by properties of the quadrature operators $Q_{I_n}$, $Q_{\bar e(t)}$ and the CFL condition \eqref{cfl} 
\begin{equation}
\label{t4} \left|\fint_{I_n}T_4(t)\; dt\right| \leq C h^{p_3+2}, \qquad  |T_5| \leq C h^{p_2+2},
\end{equation}
as the integrands are sufficiently smooth. In particular, we use the fact that the surface evolves smoothly.
Addressing the estimates for $T_1,T_2,T_3$ we will omit the time dependency as all three estimates are uniform in time.
To establish an estimate for $T_1$ we recall that we can parametrize  $e$ 
over $\bar e$ such that for the parametrisation $c_e$ inequality \eqref{le:normc} holds. We have
\begin{multline}\label{t1}
|T_1| \leq \bigg|\frac{|\bar e(t_n)|}{|e(t_n)|} -1 \bigg|\|f\|_\infty C h + \Bigg|
\int_{\bar e}\langle f^{-l}(u,x), \mu^{-l}_{K,e}(x) \rangle \left( \|c_e'(s)\| - 1\right) d\bar e\Bigg|\\
\leq\|f\|_\infty C h^{3},
\end{multline}
where $\|f\|_\infty$  denotes the supremum of $|f(u,x,t)|$ for $(x,t) \in G_T$ and $u \in \mathcal{K}.$
Next we turn to $T_3.$ Its estimate is based on the assumption that we have extended $f(u,\cdot)$
 to $\mathcal N$ smoothly 
and on the second statement of Lemma \ref{lemma:dziuk-elliott}.
This leads to
\begin{equation}\label{t3} 
|T_3| \leq \int_{\bar e} \norm{f^{-l}(u,x) - f(u,x)} \norm{\bar \mu_{\bar K,\bar e}}\, d\bar e \leq  Ch^{3}.
\end{equation}
This leaves $T_2.$ It is clear that
\begin{equation}\label{t21} |T_2| \leq Ch\max_{x \in \bar e}
\left|\langle f^{-l}(u,x), \mu_{K,e}^{-l} (x) - \bar \mu_{\bar K,\bar e}\rangle\right|.\end{equation}
Furthermore we find,
as $f$ is tangential to $\Gamma,$
\[ f^{-l}(u,x)=f_1(u,x) {\bf t}^{-l}(x) + f_2(u,x) \mu_{K,e}^{-l}(x),\]
where ${\bf t}$ is a unit tangent vector to $e$ and $f_1(u,x),f_2(u,x) \in \setR.$
Due to Lemma \ref{lem:angles1} we have
\begin{align}
\label{t22}\langle f^{-l}(u,x), \mu_{K,e}^{-l} (x)\rangle &= f_2(u,x), \\
\label{t23}\langle f^{-l}(u,x),  \bar \mu_{\bar K,\bar e}\rangle &=f_1(u,x) \landau(h^2) + f_2(u,x) + f_2(u,x) \landau(h^2).
\end{align}
Obviously it holds $|f_1(u,x)|,|f_2(u,x)| \leq \|f\|_\infty$ such that inserting \eqref{t22},\eqref{t23} into \eqref{t21} gives
\begin{equation}
\label{t2}
|T_2| \leq Ch^{3}.
\end{equation}
Now the statement of the Lemma follows from \eqref{eq:flux_est} together with \eqref{t4}, \eqref{t1}, \eqref{t3} and \eqref{t2}.
\qed\\

Our next step is to establish stability estimates for the curved and flat approximate solution.
Due to the geometry change of the surface $\Gamma$ which might act as a source term we need the following lemma.
\begin{lemma}
For every finite sequence of positive numbers $\{b_n\}_{n=1,\dots,N}$ we have
\begin{equation}
\label{eq:expest}
\prod_{n=1}^N (1 + b_n) \leq \left( 1 + \sum_{n=1}^N \frac{b_n}{N} \right)^N \leq \exp{\left(\sum_{n=1}^N b_n\right)}.
\end{equation}
\end{lemma}
\proof
From Jensen's inequality we know
\begin{equation}\label{dec1}
\sum_{n=1}^N \ln(1+ b_n) \leq N \ln \left(\sum_{n=1}^N \frac{1+b_n}{N} \right).
\end{equation}
Applying the exponential function to \eqref{dec1} gives the first inequality in \eqref{eq:expest}.
The second inequality in \eqref{eq:expest} follows from the fact that 
\[ \left( 1 + \frac{c}{N} \right)^N \leq \exp(c) \quad \forall N \in \setN , c \in \setR.\]
\qed

Now we can show a stability estimate for the curved scheme, the proof of which is mostly standard.
\begin{lemma}\label{lem:curvedstability}
Let $u_0 \in L^\infty(\Gamma(0))$. Let the
numerical flux functions of the curved scheme satisfy \eqref{def:consistency}-\eqref{def:mono}, and let the time step satisfy the CFL condition \eqref{cfl}.
Then the solution of the curved scheme fulfils 
\begin{equation}\label{dec2}
| u^{n+1}_K | \leq (1 + c |I_n|) \max\{ |u_K^n| , \max_{e \subset \partial K} \{|u_{K_e}^n|\}\} + c |I_n|
\quad \forall \, K \in \T,
\end{equation}
for some constant $c$ 
and therefore
\begin{equation} \| u^h(t) \|_{L^\infty} \leq 
( \|u_0\|_{L^\infty} +  cT) \exp(cT) \quad \forall \, 0 \leq t \leq T.\end{equation}
\end{lemma}
\proof
{Invoking the consistency of the numerical flux functions \eqref{def:consistency} we have
\[ \sum_{e \subset \partial K} |e(t_n)| f^n_{K,e}(u^n_K,u^n_K) = \fint_{I_n} \int_{K(t)} \nabla_\Gamma \cdot f(u^n_K,x,t) \, d\Gamma(t) dt.\]
Therefore, we can rewrite \eqref{eq:FVcurv} as
\begin{align*}
 u_K^{n+1} = \frac{|K(t_n)|}{|K(t_{n+1})|} & \left(  (1 - \sum_{e \subset \partial K} c_{K,e}) u_K^n + \sum_{e \subset \partial K} c_{K,e} u_{K_e}^n \right.\\
&\quad\left.-  |I_n| \fint_{I_n}\frac{1}{|K(t_{n})|}\int_{K(t)}
\nabla_{\Gamma} \cdot f(u_K^n,x,t) \, d\Gamma(t)dt \right)
\end{align*}
with
\[ c_{K,e} = \frac{|I_n|\ |e(t_n)|}{|K(t_n)|}\frac{f^n_{K,e}(u^n_K,u^n_{K_e}) - f^n_{K,e} (u^n_K,u^n_K) }{u^n_{K} -u^n_{K,e}}.\]
Due to the monotonicity of the numerical fluxes \eqref{def:mono} and the CFL condition \eqref{cfl} we have
\[ c_{K,e} \geq 0,  \sum_{e \subset \partial K} c_{K,e} \leq 1.\] 
Combining the growth condition \eqref{eq:growth} and the fact that $|K(t_n)|/|K(t_{n+1})| \leq 1+ c |I_n|$ we get 
\eqref{dec2} for another, possibly larger constant $c$.
Iteration of \eqref{dec2} implies
\begin{equation}\label{dec4} \max_{K \in \T} |u_K^n| \leq   \prod_{k=0}^{n-1} (1+ c |I_k|) \max_{K \in \T} |u_K^0|  
+ \sum_{k=0}^{n-1}  c|I_k|\prod_{j=k+1}^{n-1} (1 + c |I_j|) .  
\end{equation}
Invoking \eqref{eq:expest} we obtain from \eqref{dec4}
\begin{equation*}
 \max_{K \in \T} |u_K^n| \leq \exp(cT) \|u_0\|_{L^\infty} +  \sum_{k=0}^{n-1}  c|I_k| \exp(cT)\leq( \|u_0\|_{L^\infty} +  cT) \exp(cT).\\
\end{equation*}
\qed

As a  technical ingredient for the stability estimate of the flat scheme and the error estimate we need the following lemma whose proof is given in the appendix.
\begin{lemma}\label{lem:quotients}
 For times $t_{n},t_{n+1}$ and corresponding cells $K(t_n),$ $K(t_{n+1}),$ $\bar K(t_n),$ $\bar K(t_{n+1})$
the following estimate holds
\begin{equation}
\label{lem:quotients1} \abs{ \frac{|K(t_n)|}{|\bar K(t_{n})|} - \frac{|K(t_{n+1})|}{|\bar K(t_{n+1})|} } \leq C h |t_{n+1} - t_n|.
\end{equation}
Due to Lemma \ref{lem:eK} this implies
\begin{equation}
\label{lem:quotients2} \abs{ \frac{|\bar K(t_n)|}{|\bar K(t_{n+1})|} \frac{|K(t_{n+1})|}{|K(t_{n})|}-1 } \leq  C h |t_{n+1} - t_n|.
\end{equation}
\end{lemma}

The stability estimate for the flat scheme is a combination of the stability estimate of the curved scheme and the estimate for the difference of the fluxes.
\begin{lemma}\label{lem:flatstability}
Let $u_0 \in L^\infty(\Gamma(0))$. Let the
numerical flux functions of the curved scheme satisfy \eqref{def:consistency}-\eqref{def:mono}, and let the time step satsify the CFL condition \eqref{cfl}.
Provided \eqref{eq:fluxerror} holds for the flat numerical flux functions, then the solution of the flat scheme fulfils
\begin{equation}\label{dec5} | \bar u^{n+1}_{\bar K} | \leq (1 + 2(c+1) |I_n|) 
\max\{ |\bar u_{\bar K}^n| , \max_{\bar e \subset \partial \bar K} \{|\bar u_{\bar K_{\bar e}}^n|\}\} +  2(c+1) |I_n|  + 
 d  |I_n| h\end{equation}
for all $K \in \T$ and $0 \leq t_{n+1} \leq T$.
Here $c$ can be chosen as the same constant as in Lemma \ref{lem:curvedstability} and $d>0$ is another constant.
Therefore, for $h$ sufficiently small,
\begin{equation} \| \bar u^h(t) \|_{L^\infty} \leq (\|u_0\|_{L^\infty} + bT) \exp(b T) + 1, \quad \forall \, 0 \leq t \leq T,\end{equation}
where $b := 2(c+1).$
\end{lemma}
\proof
We have
\begin{equation}\begin{split}\label{flatstep}
 \bar u_{\bar K}^{n+1} =& \frac{|\bar K(t_n)|}{|\bar K(t_{n+1})|}  \left( \bar  u^n_K - \frac{\abs{I_n}}{|\bar K(t_{n})|}
 \sum_{\bar e \subset \partial \bar K} |\bar e(t_n)| \bar  f^n_{\bar K,\bar e}(\bar u_{\bar K}^n,\bar u_{\bar K_{\bar e}}^n) \right)\\
 =& \frac{|\bar K(t_n)|}{|\bar K(t_{n+1})|}  \left( \bar  u^n_{\bar K} - \frac{\abs{I_n}}{| K(t_{n})|}
 \sum_{e \subset \partial K} |e(t_n)| 
 f^n_{K,e}(\bar u_{\bar K}^n,\bar u_{\bar K_{\bar e}}^n) \right. \\
 &\left. + \abs{I_n} \sum_{e \subset \partial K}   \left( \frac{\abs{e(t_n)}}{\abs{K(t_n)}} f^n_{K,e}(\bar u_{\bar K}^n,\bar u_{\bar K_{\bar e}}^n)  -
  \frac{\abs{\bar e(t_n)}}{\abs{\bar K(t_n)}}  \bar f^n_{\bar K,\bar e}(\bar u_{\bar K}^n,\bar u_{\bar K_{\bar e}}^n)\right) \right).
\end{split}\end{equation}
We observe  that  because of \eqref{lem:quotients2} 
\begin{multline}
 \label{dec8}
\frac{|\bar K(t_n)|}{|\bar K(t_{n+1})|} = \frac{| K(t_n)|}{| K(t_{n+1})|} \frac{|\bar K(t_n)|}{|\bar K(t_{n+1})|} \frac{| K(t_{n+1})|}{| K(t_{n})|}\\
\leq ( 1 + c \abs{I_n} ) \cdot ( 1 + C \abs{I_n} h) \leq 1 + (c+1) \abs{I_n} ,
\end{multline}
where $c$ is the same constant as in Lemma \ref{lem:curvedstability}, for $h$ small enough.
Moreover, provided $\max_{K \in \T} |\bar u_{\bar K}^n| \leq   A+1:=(\|u_0\|_{L^\infty} + bT) \exp(b T) + 1$ we have
\begin{equation}
\frac{|\bar K(t_n)|}{|\bar K(t_{n+1})|} \left( \frac{\abs{e(t_n)}}{\abs{K(t_n)}} f^n_{K,e}(\bar u_{\bar K}^n,\bar u_{\bar K_{\bar e}}^n)  -
  \frac{\abs{\bar e(t_n)}}{\abs{\bar K(t_n)}}  \bar f^n_{\bar K,\bar e}(\bar u_{\bar K}^n,\bar u_{\bar K_{\bar e}}^n)\right)  \leq Ch
\end{equation}
because of \eqref{eq:error-e}, \eqref{eq:error-K}, and \eqref{eq:fluxerror}.
Here we have used that for $|u| ,|v| \leq A+1,$ the numerical fluxes $f^n_{K,e}(u,v)$, $\bar f^n_{\bar K,\bar e}(u,v)$ are uniformly bounded.

Provided $\max_{K \in \T}|\bar u_{\bar K}^n| \leq  A+1$ and $h, |I_n|$ sufficiently small, we 
obtain \eqref{dec5} by the same argumentation 
as in the proof of Lemma \ref{lem:curvedstability} for some $d>0$.
Note that this estimate relies on the growth condition \eqref{eq:growth}.
As obviously $\| \bar u^h(0) \|_{L^\infty} \leq A + 1$ we have by induction
\begin{multline}\label{dec6} \max_{K \in \T} |\bar u_K^n| \leq   \prod_{k=0}^{n-1} (1+ b |I_k|) \max_{K \in \T} |\bar u_K^0|  
+ \sum_{k=0}^{n-1}  (b|I_k| + d h |I_k| )\prod_{j=k+1}^{n-1} (1 + b |I_j|)   \\
\leq \underbrace{( \| u_0\|_{L^\infty} + bT) \exp(bT)}_{\leq A} + d T \exp(bT) h,
\end{multline}
where $b=2(c+1)$.
Equation \eqref{dec6}
shows that our induction hypothesis, $\max_{K \in \T}|\bar u_{\bar K}^n| \leq  A+1$,
also holds for the next time step provided $ h < \frac{1}{ \exp(bT)d T}$ and $t_n\leq T.$
This implies that \eqref{dec5} and \eqref{dec6}  in fact hold for all $t_n \leq T$.
Thus, provided $h$ is small enough, the assertion of the lemma follows by induction.
\qed

In addition we need the fact that the curved scheme satisfies a discrete $L^1$-contraction property.
\begin{lemma}\label{lem:contraction}
 For given data $u_K^n$ and $v_K^n,$ let $u_K^{n+1}$ and $v_K^{n+1},$  be defined according to 
 the curved finite volume scheme 
\eqref{eq:FVcurv}.
If the corresponding numerical flux functions and the time step satisfy \eqref{def:consistency}-\eqref{cfl}, then 
\[ \sum_K |K(t_{n+1})| | u_K^{n+1} - v_K^{n+1} | \leq \sum_K |K(t_{n})| | u_K^{n} - v_K^{n} |.\]
\end{lemma}
As the appropriate cell weights (which depend on the time step) appear in the scheme as well as in the $L^1$-norm the
 proof is analogous to the proof of the discrete $L^1$-contraction property for finite volume schemes in Euclidean
space, cf. \cite{CCL94}.

For the difference between the curved and flat approximate solutions we obtain the following estimate.
\begin{theorem}
\label{thm:main}
For initial data $u_0 \in L^\infty(\Gamma_h(0))$, let $u^h$ denote the solution of the curved finite volume scheme \eqref{eq:FVcurv}
and let $\bar u^h$ denote the solution of the flat finite volume scheme \eqref{eq:FVflat}.
Let, in addition, the quadrature operators $Q_{\bar K(0)}$ and the initial data $u_0$ be such that
\begin{align}\label{eq:mainthmass}
 \|u^h(0)-\bar u^h(0)\|_{L^1(\Gamma(0))} \leq C \,h
\end{align}
for some constant $C$.
If the curved numerical flux functions and the time step satisfy \eqref{def:consistency}-\eqref{cfl}, and additionally,
\eqref{eq:fluxerror} holds for the flat numerical flux functions,
 then, for fixed $T>0$,
the difference between $u^h$
 and $\bar u^h$ satisfies
 \begin{equation}\label{eq:errorestimate}
 \norm{ u^h(T) - \bar{u}^h(T) }_{L^1(\Gamma(T))} \leq  C \, h,
 \end{equation}
for some constant $C$ depending on $G_T,f,u_0$.
\end{theorem}

\begin{remark}\label{rem:error}
The curved approximate solution converges to the entropy solution of \eqref{eq:ConLaw}-\eqref{eq:IVConLaw}
with a convergence rate of $\mathcal O(h^{1/4})$, cf. \cite{GW12}.
 Hence, invoking Theorem \ref{thm:main} the same kind of error bound holds for the flat approximate solution.
\end{remark}
As the Lax-Friedrichs numerical flux functions from Section \ref{sec:fvs} satisfy assumptions
\eqref{def:consistency}-\eqref{def:mono} and \eqref{eq:fluxerror} due to Lemma \ref{lem:flux_est}
we obtain the following corollary.
\begin{corollary}
Let the numerical solutions $u^h$ and $\bar u^h$ be defined with the Lax-Friedrichs numerical fluxes 
\eqref{curv_num_flux} and \eqref{eq:LFflat}. Let the time step sizes satisfy the CFL condition \eqref{cfl}. Then, the error estimate \eqref{eq:errorestimate}
holds.
\end{corollary}

\proof[of Theorem \ref{thm:main}] 
Let $n\in\setN$ be such that $T \in [t_n, t_{n+1})$, then we have
\begin{align*}
\| u^h(T) - &\bar{u}^h(T) \|_{L^1(\Gamma)} =  \sum_K\abs{K(t_{n+1})} \abs{ u_{K}^{n+1} - \bar u _{\bar K}^{n+1} } \\
=   \sum_K \Bigg|& |K(t_n)| u_{K}^{n}  - |I_n| \sum_{e\subset \partial K}  
 \abs{e(t_n)}  f^n_{K,e} ( u_{K}^{n} ,u_{K_e}^{n}) \\
	&- \frac{|K(t_{n+1})|\, | \bar K(t_n)|}{|\bar K(t_{n+1})|}\bar u _{\bar K}^{n}
 + |I_n|
\frac{|K(t_{n+1})|}{|\bar K(t_{n+1})|}\sum_{e \subset \partial K} \abs{\bar e(t_n)}  \bar f^n_{\bar K,\bar e}
 ( \bar u_{\bar K}^{n} ,\bar u_{\bar K_{\bar e}}^{n}) \Bigg|  \\
\leq R_1 &+ R_2 + R_3 + R_4 + R_5,
\end{align*}
where
\begin{align*}
R_1:= & \sum_K \Big| |K(t_n)| u_{K}^{n}  - |I_n| \sum_{e\subset \partial K}  
 \abs{e(t_n)}  f^n_{K,e} ( u_{K}^{n} ,u_{K_e}^{n})\\
 & \qquad-|K(t_n)| \bar u_{\bar K}^{n}  + |I_n| \sum_{e\subset \partial K}  
 \abs{e(t_n)}  f^n_{K,e} ( \bar u_{\bar K}^{n} ,\bar u_{\bar K_{\bar e}}^{n}) \Big|\\
R_2:=&\sum_K \Big||K(t_n)| \bar u_{\bar K}^n -  
 \frac{|K(t_{n+1})|\, | \bar K(t_n)|}{|\bar K(t_{n+1})|}\bar u_{\bar K}^n  \Big| \\
R_3:=&\sum_K  |I_n| \sum_{e\subset \partial K}  \abs{ \abs{e(t_n)} - \abs{\bar e(t_n)} }
 \abs{f^n_{K,e} (\bar u_{\bar K}^{n} ,\bar  u_{{\bar K}_{\bar e}}^{n})}\\
R_4:=&\sum_K\left| \left( 1 - \frac{|K(t_{n+1})|}{|\bar K(t_{n+1})|} \right) |I_n| \sum_{e\subset \partial K}  \abs{\bar e(t_n)} 
 \abs{f^n_{K,e} (\bar u_{\bar K}^{n} ,\bar  u_{{\bar K}_{\bar e}}^{n})} \right|\\
R_5:=&\sum_K   |I_n|   \frac{|K(t_{n+1})|}{|\bar K(t_{n+1})|}    \sum_{e\subset \partial K}  
 \abs{\bar e(t_n)} \abs{  f^n_{ K, e} (\bar u_{\bar K}^{n} ,\bar  u_{{\bar K}_{\bar e}}^{n}) 
	-  \bar f^n_{\bar K,\bar e} (\bar u_{\bar K}^{n} ,\bar  u_{{\bar K}_{\bar e}}^{n})   }.
\end{align*}
According to Lemma \ref{lem:contraction} the curved finite volume scheme satisfies the $L^1$-contraction property and therefore
\begin{align*}
 R_1 \leq  \sum_{K} \abs{K(t_n)} \abs{ u_{K}^{n} - \bar u _{\bar K}^{n} }.
\end{align*}
The term $R_2$ can be estimated using \eqref{lem:quotients1}, we get
\begin{equation}
 R_2 \leq \sum_K \abs{ \bar u_{\bar K}^n} |\bar K(t_n)| \underbrace{ \left| \frac{|K(t_n)|}{|\bar K(t_n)|}-     
\frac{|K(t_{n+1})|}{|\bar K(t_{n+1})|}\right|}_{\leq C |I_n| h}   \leq C |I_n| h.
\end{equation}
Applying Lemma \ref{lem:eK} and assumption \eqref{ineq:angle-assumption} together with Remark \ref{rem:angles} we get 
\begin{equation}
 R_3,R_4 \leq \sum_K |I_n| \sum_{e \subset \partial K} C h^3  
\abs{f^n_{K,e} (\bar u_{\bar K}^{n} ,\bar  u_{{\bar K}_{\bar e}}^{n})}
\leq C |I_n| h.\\
\end{equation}
Based on Lemma \ref{lem:eK}, assumption \eqref{ineq:angle-assumption}, Remark \ref{rem:angles}
and \eqref{eq:fluxerror}
we have
\begin{equation}
R_5 \leq C \sum_K |I_n| \sum_{e \subset \partial K} h^3 \leq C |I_n| h.
\end{equation}
%
Combining these estimates we thus obtain by iteration
\begin{align*}
 \norm{ u^h(T) - \bar{u}^h(T) }_{L^1(\Gamma)} &=  \sum_K\abs{K(t_{n+1})} \abs{ u_{K}^{n+1} - \bar u _{\bar K}^{n+1} }  \\
&\leq  \sum_K\abs{K(t_n)} \abs{ u_{K}^{n} - \bar u _{\bar K}^{n} } + C |I_n|\; h\\
&\leq  \sum_K\abs{K(0)} \abs{ u_{K}^{0} - \bar u _{\bar K}^{0} } + C T\; h\\
&\leq C (T+1) h,
\end{align*}
where the last step follows with \eqref{eq:mainthmass}.
\qed

\section{Numerical Experiments}
\label{sec:numerics}

Numerical investigations based on the finite volume schemes defined in Section \ref{sec:fvs} are presented in this section.
The upshot of our experiments is three-fold. Firstly, under the present assumptions the order of convergence stated in Theorem \ref{thm:main}
is optimal. This is demonstrated by Test Problem \ref{tp:independent}.
Secondly, all of our experiments which include a sufficiently large numerical viscosity, 
i.e. $\lambda\in\operatorname\Theta (1)$ in \eqref{curv_num_flux},
lead to a considerably higher experimental order of convergence (EOC) between 1 and 2 for the $L^1$-difference between the 
flat and the curved approximate solution. Thirdly, the application of a finite volume scheme of second order 
to Test Problem \ref{tp:independent} demonstrates that orders of convergence higher than 1 are not to be
expected in general, if the geometry is not approximated sufficiently well, see Test Problem \ref{tp:hofvtp1}.
In the following we will present several test cases. Thereafter, we will mention some implementation aspects.

\subsection{Test Problems}

All test cases except Test Problems \ref{tp:shrinkingsphere} and \ref{tp:movingsurface} use the geometrical setting
$G_T=\mathbb S^2\times [0,1]$, i.e. $\Gamma(t)=\mathbb S^2$ for all $t\in [0,T]$, and $T=1$. This is due to the fact,
that we are able to compute the exact curved quantities only in this or similarly simple settings. In addition,
let us fix the vector fields
$V(x)=\frac{2\pi}{\|x\|}(x_2,-x_1,0)^T$
and $W(x)=\frac{2\pi}{\|x\|}(-x_3,0,x_1)^T$ for
$x\in \setR^3\backslash \{0\}$.

\begin{testproblem}[$u$-independent flux function]
We choose $f=V$ as the flux function. Since $f$
neither depends on $t$ nor on $u$ and is divergence-free on $\mathbb S^2$ any initial datum $u_0: \mathbb S^2\rightarrow \mathbb R$ is
a stationary solution of the corresponding initial value problem \eqref{eq:ConLaw}-\eqref{eq:IVConLaw}.
For initial values identically to zero the curved scheme conserves this stationary solution.
Thus, the error between the curved and the flat approximate solution
is equal to the error between the flat approximate solution and the exact solution.
The results for this test case for $\lambda=0$ are plotted in Table \ref{tab:1}.
Note that due to $\partial_u f=0$ the numerical flux functions are monotone.
This experiment shows, that under the assumptions from our convergence analysis $\mathcal O(h)$ is indeed the optimal order of convergence.

However, if we modify the numerical diffusion by setting $\lambda=\pi$ in the numerical flux functions we achieve EOCs between $1$ and $2$
as can be seen in Table \ref{tab:1}, as well.
\label{tp:independent}
\end{testproblem}

\begin{table}
\centering\scalebox{0.95}{
\begin{tabular}{|c|@{}p{0.2em}@{}|c|c|@{}p{0.2em}@{}|c|c|}
\cline{3-4}\cline{6-7}
 \multicolumn{1}{c}{} &  & \multicolumn{2}{|c|}{Test Problem \ref{tp:independent}, $\lambda = 0$} &  & \multicolumn{2}{|c|}{Test Problem \ref{tp:independent}, $\lambda = \pi$}\\ 
\cline{1-2}\cline{3-4}\cline{5-7}
 \multicolumn{1}{|c|}{level} &  & \multicolumn{1}{|c|}{$L^1$-difference} & \multicolumn{1}{|c|}{EOC} &  & \multicolumn{1}{|c|}{$L^1$-difference} & \multicolumn{1}{|c|}{EOC}\\
\cline{1-2}\cline{3-4}\cline{5-7}
$0$ & & $0.758314$ & ---     &      & $0.0119577$ & --- \\
$1$ & & $0.437173$ & $0.805$ &      & $0.0050082$ & $1.271$ \\
$2$ & & $0.231999$ & $0.917$ &      & $0.0020877$ & $1.266$ \\
$3$ & & $0.119190$ & $0.962$ &      & $0.0008286$ & $1.334$ \\
$4$ & & $0.060372$ & $0.982$ &      & $0.0003137$ & $1.402$ \\
$5$ & & $0.030378$ & $0.991$ &      & $0.0001165$ & $1.429$ \\
$6$ & & $0.015237$ & $0.995$ &      & $0.0000439$ & $1.408$ \\
\cline{1-2}\cline{3-4}\cline{5-7}
\end{tabular}}
\caption{$L^1$-difference and EOCs between curved approximate solution $u^h(T)$ and
flat approximate solution $\bar u^h(T)$ from Test Problem \ref{tp:independent} for
different values $\lambda$ of numerical diffusion.}
\label{tab:1}
\end{table}

\begin{testproblem}[Advection across the poles]
\label{tp:linear}
Let the flux function $f$ be defined by $f(u,x)= u W(x)$ for $x\in \mathbb S^2$. Initial values
are given by $u_0(x)= \mathbbmss 1_{\{x_1>0.15\}}(x)$. In order to get monotone numerical flux functions
we set $\lambda=\frac{1}{2}\|\partial_u f \|_\infty=\pi$. For this test case we obtain EOCs of almost $2$, cf. Table \ref{tab:2}.
\end{testproblem}

\begin{testproblem}[Burgers along the latitudes]
\label{tp:Burgers}
We choose a flux function of Burgers-type $f=f(u,x)=1/2 u^2 V(x)$ for $x\in \mathbb S^2$
and initial values $u_0(x)= \mathbbmss 1_{\{x_1>0.15\}}(x)$. In order to get monotone numerical flux functions
we set $\lambda=\frac{1}{2}\|\partial_u f \|_\infty=\pi$ and obtain EOCs of 
almost $2$, cf. Table \ref{tab:2}. 
\end{testproblem}

\begin{testproblem}[Fully two-dimensional problem]
\label{tp:twodim}
In this test problem we consider a flux function $f$ such that the corresponding initial value problem
is not equivalent to a family of one-dimensional problems. Note that the flux functions from the previous test problems
have been of one-dimensional nature. To this end we define $f(u,x)=u V(x) + 1/2u^2 W(x)$ for $x\in\mathbb S^2$ with
initial values $u_0(x)= \mathbbmss 1_{\{x_1>0.15\}}(x)$ and observe EOCs of almost $2$, cf. Table \ref{tab:2}.
\end{testproblem}

\begin{table}
\centering\scalebox{0.95}{
\begin{tabular}{|c|@{}p{0.2em}@{}|c|c|@{}p{0.2em}@{}|c|c|@{}p{0.2em}@{}|c|c|}
\cline{3-4}\cline{6-7}\cline{9-10}
 \multicolumn{1}{c}{} &  & \multicolumn{2}{|c|}{Test Problem \ref{tp:linear}} &  & \multicolumn{2}{|c|}{Test Problem \ref{tp:Burgers}}&  & \multicolumn{2}{|c|}{Test Problem \ref{tp:twodim}}\\ 
\cline{1-2}\cline{3-4}\cline{5-7}\cline{8-10}
 \multicolumn{1}{|c|}{level} &  & \multicolumn{1}{|c|}{$L^1$-difference} & \multicolumn{1}{|c|}{EOC} &  & \multicolumn{1}{|c|}{$L^1$-difference} & \multicolumn{1}{|c|}{EOC}&  & \multicolumn{1}{|c|}{$L^1$-difference} & \multicolumn{1}{|c|}{EOC}\\
\cline{1-2}\cline{3-4}\cline{5-7}\cline{8-10}
$0$ & & $0.112518$ & --- &      & $0.0370831$ & --- &      & $0.115867$ & --- \\
$1$ & & $0.039167$ & $1.541$ &      & $0.0133379$ & $1.494$ &      & $0.035202$ & $1.740$ \\
$2$ & & $0.011223$ & $1.809$ &      & $0.0040350$ & $1.730$ &      & $0.009566$ & $1.886$ \\
$3$ & & $0.002984$ & $1.913$ &      & $0.0011216$ & $1.848$ &      & $0.002475$ & $1.952$ \\
$4$ & & $0.000772$ & $1.951$ &      & $0.0002992$ & $1.907$ &      & $0.000630$ & $1.974$ \\
$5$ & & $0.000197$ & $1.970$ &      & $0.0000778$ & $1.943$ &      & $0.000159$ & $1.986$ \\
$6$ & & $0.000053$ & $1.894$ &      & $0.0000199$ & $1.967$ &      & $0.000040$ & $1.991$ \\
\cline{1-2}\cline{3-4}\cline{5-7}\cline{8-10}
\end{tabular}}
\caption{$L^1$-difference and EOCs between curved approximate solution $u^h(T)$ and
flat approximate solution $\bar u^h(T)$ from Test Problems \ref{tp:linear}, \ref{tp:Burgers} and \ref{tp:twodim}.}
\label{tab:2}
\end{table}

\begin{testproblem}[
2nd order scheme applied to Test Problem \ref{tp:independent}]
\label{tp:hofvtp1}
The\\ motivation of this test problem is to show that in general
even higher order schemes, which are based on the flat finite volume schemes,
are not able to achieve higher order convergence rates for smooth data.
To this end, we apply a second order finite volume scheme (which is validated in Test Problem \ref{tp:hofvtpval})
to Test Problem \ref{tp:independent}. This scheme is based on the flat finite volume scheme of first order (cf. Subsection \ref{subsec:fv_flat})
with $\lambda=0$
enhanced with a linear reconstruction and a second order Runge-Kutta method for time evolution. In Table
\ref{tab:3} we observe EOCs of almost $1$.
Indeed, the application of a second order finite
volume scheme to Test Problem \ref{tp:independent} gives almost exactly the same convergence rates as a first order scheme since the 
linear reconstruction on each cell does not affect the numerical flux functions as $f$ is independent of $u$.
Note that both schemes are not identical, e.g. the time integrations are different, explaining the slight deviation of the EOCs.
We like to point out that we do not have to compute the curved approximate solution as it coincides with the 
(constant) exact solution.
\end{testproblem}

\begin{testproblem}[Validation of the 
2nd order scheme]
\label{tp:hofvtpval}
 This test problem serves as validation of the second order finite volume scheme. 
 We consider smooth initial values 
 \begin{align*}
u_0(x):=\frac{1}{10} \ \mathbbmss 1_{\{r(x)<1\}}(x)\ \exp \left( \frac{-2 \left(  1+r^2(x)\right)}{\left(1-r^2(x)\right)^2}   \right)  
 \end{align*}
 with $r(x):=\frac{|x_0-x|}{0.74}$ and $x_0:=(1,0,0)^T$
 and a flux function $f(u,x):= u V(x)$, which 
 transports the initial values around the sphere.
 For the error between the flat second order finite volume scheme (see Test Problem \ref{tp:hofvtp1}) and the exact solution 
 EOCs significantly higher than $1$ are shown in Table \ref{tab:3}.
\end{testproblem}

\begin{table}
\centering\scalebox{0.95}{
\begin{tabular}{|c|@{}p{0.2em}@{}|c|c|@{}p{0.2em}@{}|c|c|@{}p{0.2em}@{}|c|c|}
\cline{3-4}\cline{6-7}\cline{9-10}
 \multicolumn{1}{c}{} &  & \multicolumn{2}{|c|}{Test Problem \ref{tp:hofvtp1}} &  & \multicolumn{2}{|c|}{Test Problem \ref{tp:hofvtpval}}&  & \multicolumn{2}{|c|}{Test Problem \ref{tp:shrinkingsphere}}\\ 
\cline{1-2}\cline{3-4}\cline{5-7}\cline{8-10}
 \multicolumn{1}{|c|}{level} &  & \multicolumn{1}{|c|}{$L^1$-difference} & \multicolumn{1}{|c|}{EOC} &  & \multicolumn{1}{|c|}{$L^1$-error} & \multicolumn{1}{|c|}{EOC}&  & \multicolumn{1}{|c|}{$L^1$-error} & \multicolumn{1}{|c|}{EOC}\\
\cline{1-2}\cline{3-4}\cline{5-7}\cline{8-10}
$0$ & & $0.777427$ & ---     &      & $0.00362492$ & ---     &      & $2.54384$ & --- \\
$1$ & & $0.444068$ & $0.818$ &      & $0.00243102$ & $0.584$ &      & $1.99697$ & $0.354$ \\
$2$ & & $0.233521$ & $0.930$ &      & $0.00112593$ & $1.114$ &      & $1.49537$ & $0.419$ \\
$3$ & & $0.119553$ & $0.967$ &      & $0.00034768$ & $1.697$ &      & $1.09177$ & $0.454$ \\
$4$ & & $0.060461$ & $0.984$ &      & $0.00012006$ & $1.534$ &      & $0.78774$ & $0.471$ \\
$5$ & & $0.030400$ & $0.992$ &      & $0.00003713$ & $1.693$ &      & $0.56469$ & $0.480$ \\
$6$ & & $0.015242$ & $0.996$ &      & $0.00001116$ & $1.734$ &      & $0.40312$ & $0.486$ \\
\cline{1-2}\cline{3-4}\cline{5-7}\cline{8-10}
\end{tabular}}
\caption{$L^1$-difference and EOCs between a second order curved approximate solution (which equals 
the exact solution in this case) and a second order flat approximate solution from Test Problem \ref{tp:hofvtp1},
$L^1$-error and EOCs between the exact solution from Test Problem \ref{tp:hofvtpval} and its approximation by the second order finite volume scheme
and $L^1$-error and EOCs between the exact solution from Test Problem \ref{tp:shrinkingsphere} and its approximation by the
flat approximate solution using the Lax-Friedrichs numerical flux functions from \eqref{eq:LFflat}.}
\label{tab:3}
\end{table}

\begin{testproblem}[Shrinking Sphere]\label{tp:shrinkingsphere}
In order to validate the convergence rate of the error between the flat approximate solution and
 the entropy solution from Remark \ref{rem:error}
 we consider a linear transport problem on a shrinking sphere that was introduced in \cite{DKM13}.
 In the conservation law \eqref{eq:ConLaw} let 
 \begin{align*}
 \Gamma(t):= \exp(-t)\ \mathbb S^2 \ \text{ and }\ f(u,x,t):=-u\ V(x)
 \end{align*}
for $t\in [0,1]$ and $x\in \Gamma(t)$. Analogously to \cite{DKM13} one sees that the function $u$, expressed in spherical coordinates 
$(\varphi,\theta)\in \mathbb (0,2\pi)\times (0,\pi)$  by
\begin{align*}
u(\varphi,\theta,t) &:= \exp(2t) \widetilde u(\varphi -  2\pi (\exp(t)-1)) \widehat u(\theta),\\
\widetilde u(\varphi)&:=\mathbbmss 1_{\left\{\varphi <\pi/2\right\} }(\varphi),\\
  \widehat u (\theta) & := \mathbbmss 1_{\left\{|\theta-\pi/2|<\pi/4\right\}}(\theta)
\end{align*}
solves \eqref{eq:ConLaw} for initial values $u_0(\varphi,\theta) := \widetilde u(\varphi) \widehat u(\theta)$.
For the error between the exact solution and the flat approximate solution at end time $T=1$ we observe, identically to  
similar problems in the Euclidean space, EOCs of almost $0.5$, cf. Table \ref{tab:3}.
\end{testproblem}

\begin{testproblem}[Deforming Torus]
 \label{tp:movingsurface}
We consider a deforming torus as computational domain $\Gamma$ and $T=4$ as final time.
Within the time interval $[0,2]$ the right half of the torus undergoes compression whereas the left half is stretched, while
$\Gamma(t)$ remains constant for $t\in[2,4]$.
We choose a Burgers-type flux function $f=f(u,x)=\frac{1}{2} u^2 (x_2,-x_1,0)^T$ and
constant initial values $u_0\equiv 1$. The time step size is chosen dynamically for each
time step such that stability is guaranteed.
In Figure \ref{fig:movingsurface} the numerical solution is shown at four different
times. Note that in spite of the constant initial values, a shock wave
is induced due to the change of geometry (compression and rarefaction) and
the nonlinearity of the flux function.
\end{testproblem}

\begin{figure}[htbp]
  \centering
  \subfigure[$t=0$.]{
    \includegraphics[width=0.45\textwidth]{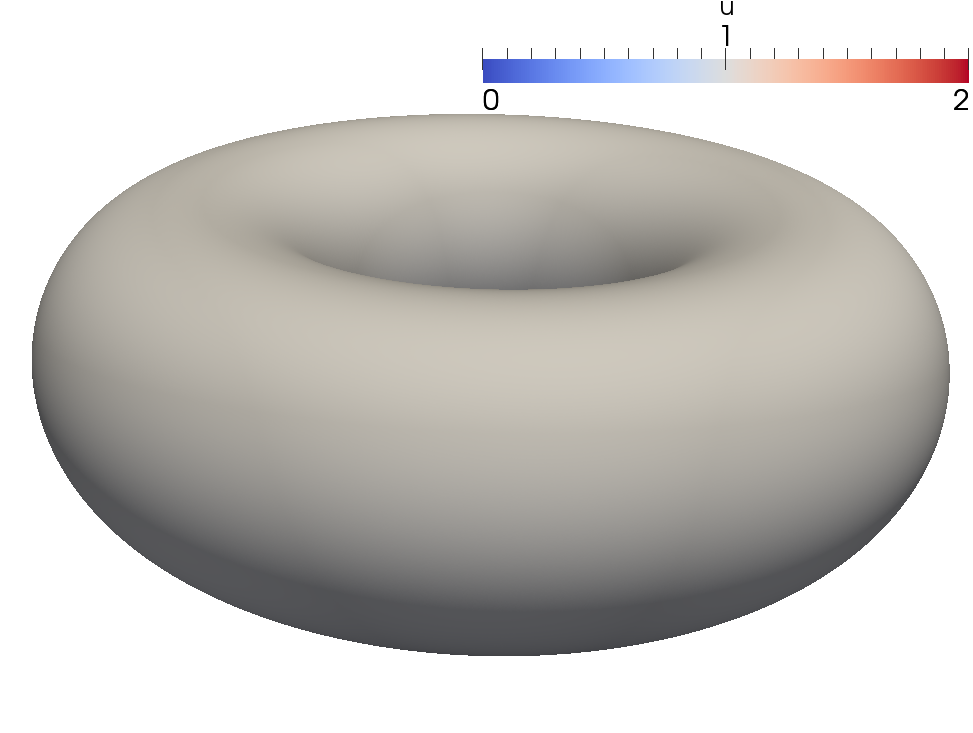} 
  }
  \subfigure[$t=1.07$.]{
    \includegraphics[width=0.45\textwidth]{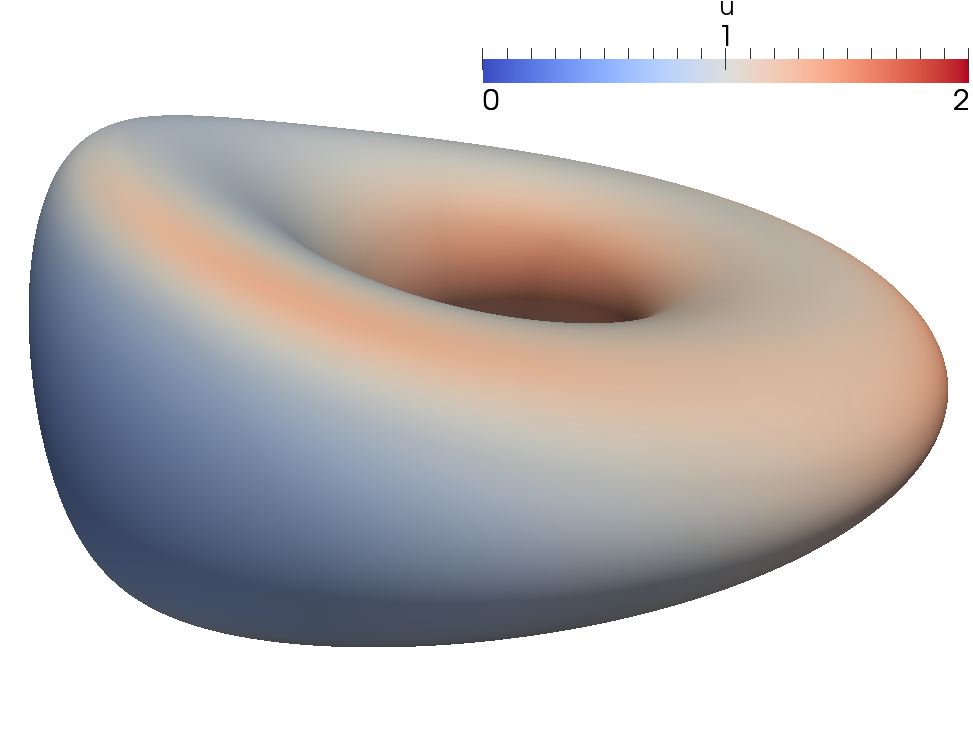} 
  }
  \subfigure[$t=2.36$.]{
    \includegraphics[width=0.45\textwidth]{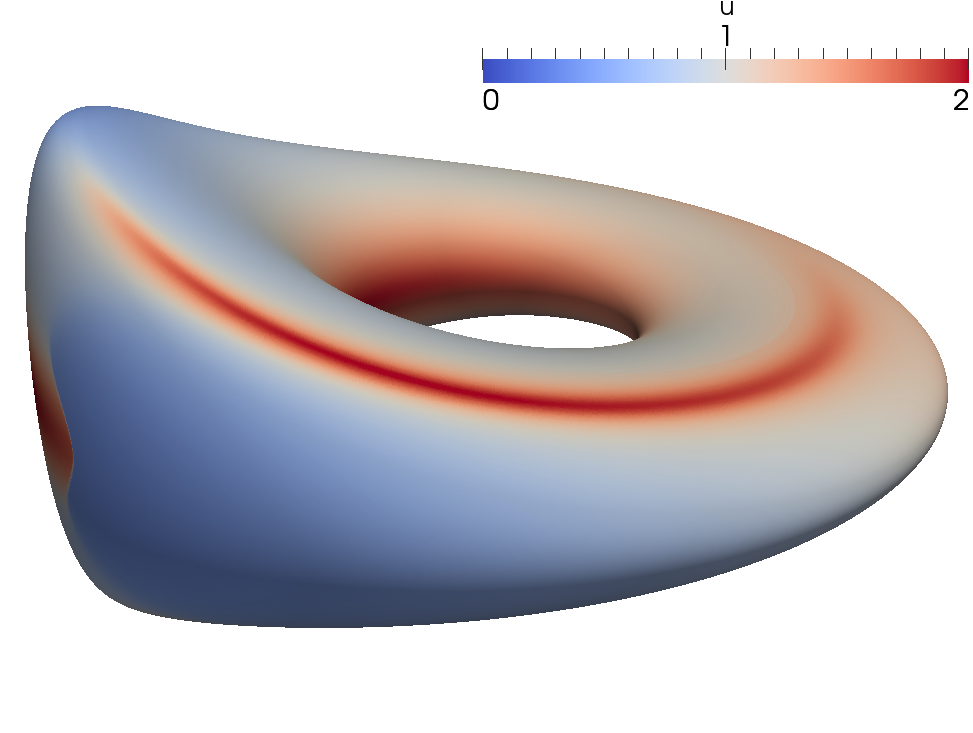} 
  }
  \subfigure[$t=4$.]{
   \includegraphics[width=0.45\textwidth]{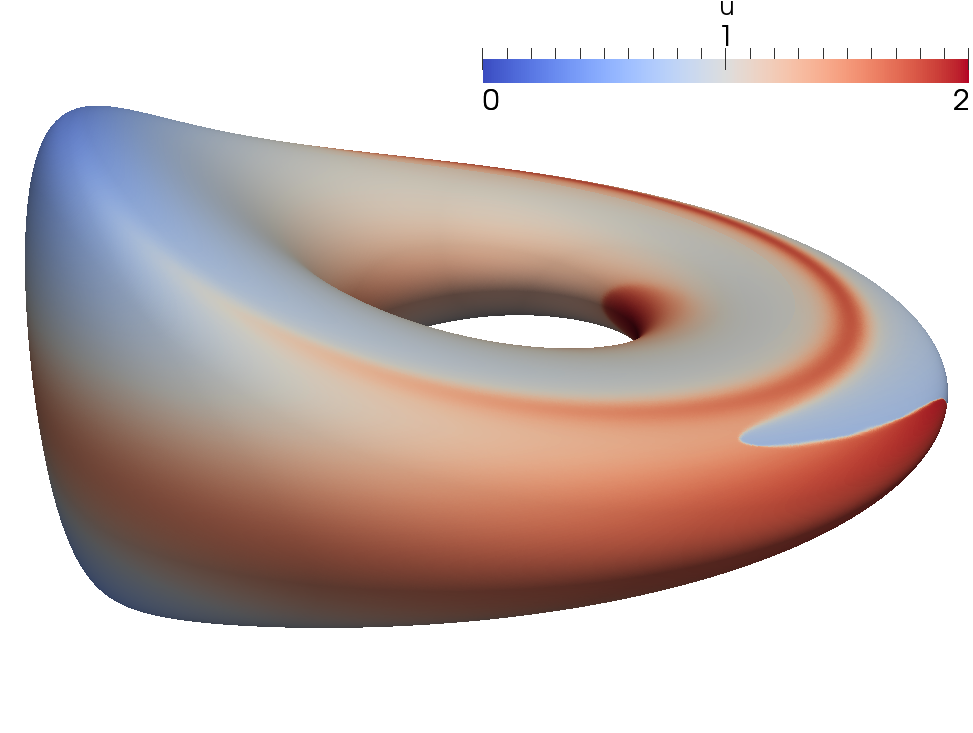} 
  }
  \caption{Flat approximate solution for Test Problem \ref{tp:movingsurface} for four different times.
  The computation was performed on a deforming polyhedron consisting of about 3 million triangles.}
  \label{fig:movingsurface}
\end{figure}

\subsection{Implementation Aspects}

\subsubsection{Software}
All simulations have been performed within the DUNE-FEM module (see \cite{DKNO10} and the references therein) 
which is based on the Distributed 
and Unified Numerics Environment (DUNE) using ALUGRID \cite{KN} as grid implementation.
The figures have been created with ParaView. 
As coarsest grid approximating the sphere we use an unstructured grid consisting of $632$ triangles, see Figure \ref{fig:sphere0}.
For finer computations we refine the coarse macro grid (level $0$)
and obtain up to $2.5$ million triangles for the finest grid (level $6$) whose vertices are projected onto the sphere, cf. Table \ref{tab:sphere}.

\begin{figure}[htbp]
\begin{minipage}[b]{0.47\textwidth}
\centering\scalebox{0.95}{
\begin{tabular}{|c|r|r|}
\cline{1-3}
\multicolumn{1}{|c|}{level} & \multicolumn{1}{|c|}{h} & \multicolumn{1}{|c|}{size}  \\
\cline{1-3}
$0$ & $0.311151$ & $632$  \\
$1$ & $0.156914$ & $2528$ \\
$2$ & $0.078628$ & $10112$ \\
$3$ & $0.039335$ & $40448$ \\
$4$ & $0.019670$ & $161792$ \\
$5$ & $0.009835$ & $647168$ \\
$6$ & $0.004918$ & $2588672$ \\
\cline{1-3}
\end{tabular}}
\captionof{table}{Different refinement levels of the sphere grid.}
\label{tab:sphere}
\end{minipage}
\begin{minipage}[b]{0.47\textwidth}
\centering
\includegraphics[width=0.95\textwidth]{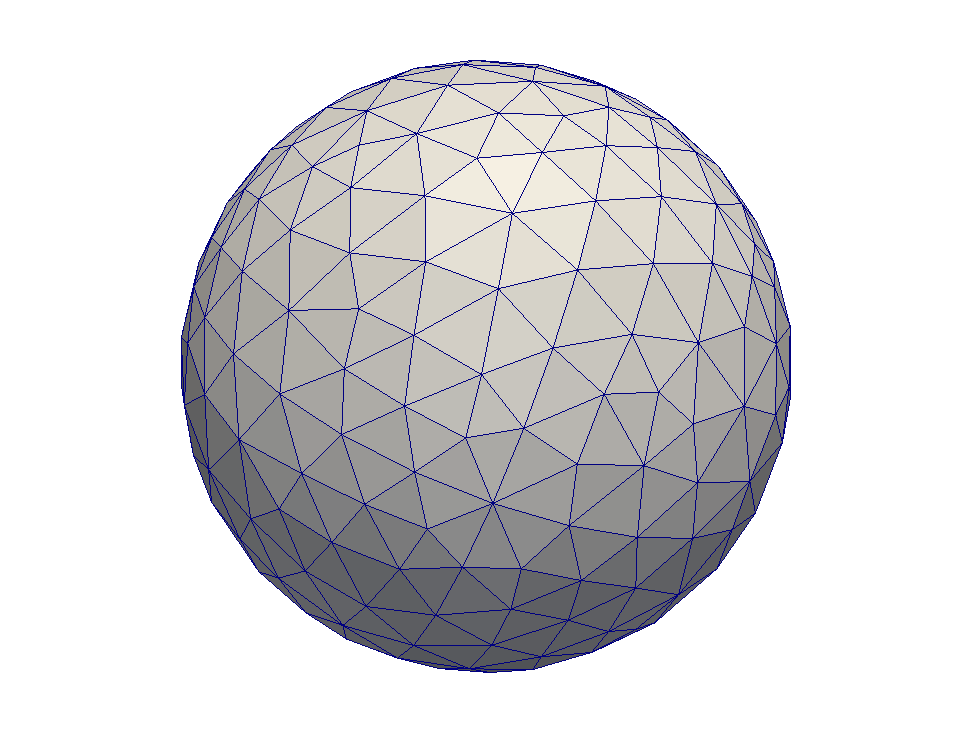}
\caption{The sphere grid of level $0$.}
\label{fig:sphere0}
\end{minipage}
\end{figure}

\subsubsection{Exact Computation of Spherical Volume}
For the curved finite volume scheme on the sphere the exact outer conormals, exact lengths of boundary segments 
and exact volumes of 
spherical triangles need to be computed. While the computation of the former two quantities is an easy geometric exercise,
we use the formula from \cite{OS83} for the computation of the latter.

\subsubsection{Exact Computation of Numerical Flux Functions}

For the exact evaluation of the numerical flux function corresponding to an edge $e$
of a grid cell $K$, quantities of the form
 $\fint_{e} \langle  V, \mu_{ K, e} \rangle\, de$
have to be computed.
Note that $V$ can be written as $V = \nu \times \nabla h_V$ with $h_V(x)=2\pi x_3$,
where $\nu(x):=x$ denotes the outer unit normal to $\mathbb S^2$. As a result, similar to \cite{DKM13},
we deduce
\begin{align*}
 \fint_{e} \langle  V, \mu_{ K, e} \rangle\, de = 
 \fint_{e} \langle  \mu_{ K, e} \times \nu, \nabla h_V \rangle\, de.
\end{align*}
As $\mu_{ K, e} \times \nu$ is a unit tangent vector to $e$, the integrand is a directional derivative along $e$
and thus the integral can be computed by the evaluation at the endpoints of $e$.
Obviously, the same applies to $W$ with $h_W(x)=2\pi x_2$ and 
$W = \nu \times \nabla h_W$.

\subsubsection{Computation of $L^1$-Norms}
We remark that the $L^1$-differences between the flat and the curved approximate solutions are computed 
on the triangulation $\Gamma_h$. This does not have any influence on the convergence rates.

\begin{acknowledgements}
We gratefully acknowledge that the work of Thomas M\"uller was supported by the 
German Research Foundation (DFG) via SFB TR 71 `Geometric Partial
Differential Equations' and by the German National Academic Foundation (Studienstiftung des Deutschen Volkes). 
Jan Giesselmann would like to thank the German Research Foundation (DFG) for financial support of the project
`Modeling and sharp interface limits of local and non-local generalized
Navier--Stokes--Korteweg Systems'.

The authors would like to express their gratitude to the two
anomymous referees for their constructive suggestions to improve
this work.
\end{acknowledgements}

\section*{Appendix}
Here we give the proof of Lemma \ref{lem:quotients}.
\proof
It is sufficient to show
\begin{equation}\label{eq:db}
 \Big| \frac{d}{dt} \frac{|K(t)|}{ |\bar K(t)|} \Big| \leq Ch.
\end{equation}
We have
\begin{multline}
 \frac{d}{dt} \frac{|K(t)|}{|\bar  K(t)|} =  \frac{d}{dt} \frac{|K(t)|- |\bar K(t)| }{|\bar K(t)|}  \\
= \frac{(|K(t)|- |\bar K(t)| )_t } {|\bar K(t)|}  + \frac{1}{ |\bar K(t)|} \Big( \frac{|K(t)|}{ |\bar K(t)|} -1 \Big)  |\bar K(t)|_t
\end{multline}
so that, due to Lemma \ref{lem:eK}, it suffices to show 
\begin{equation}\label{eq:app1}
 \Big|\frac{(|K(t)|- |\bar K(t)| )_t } { |\bar K(t)|} \Big|\leq C h \quad \text{ and } \quad \big| |\bar K(t)|_t\big| \leq Ch.
\end{equation}
The estimate \eqref{eq:app1}$_2$ is immediate, so we turn our attention to \eqref{eq:app1}$_1$.
Let us assume $t \in [t_n,t_{n+1}]$ and let $\bar K(t_n)$ be the convex hull of the vertices $v_0:=(0,0,0),\, v_1:= (h,0,0),\,v_2:= (x,y,0).$ We define
\[ \Phi^n(\cdot, t) : \Gamma(t_n) \rightarrow \Gamma(t), \quad 
 \Phi^n(\cdot,t):= \Phi(\cdot,t) \circ \Phi(\cdot,t_n)^{-1},
\]
such that $\Phi^n(\cdot,t_n)$ is the identity map.
We denote the canonical projection $\bar K(t_n) \to K(t_n)$ by $c$ and abbreviate $\Phi^n \circ c$ by $\Phi^n_c$.
The scaled directional derivatives are denoted $\partial_{v_1} := h \partial_{x_1}$ and $\partial_{v_2}:= x\partial_{x_1} + y\partial_{x_2}.$ Then,
\begin{multline}
 (|K(t)|- |\bar K(t)| )_t  = 
\frac{d}{dt} \Big( \frac{1}{hy} \int_{\bar K(t_n)} \| \partial_{v_1} \Phi^n_c(x,t) \times \partial_{v_2} \Phi^n_c(x,t)\| \, d\bar K(t_n) \\
-\frac{1}{2} \| (\Phi^n(v_1,t) - \Phi^n(v_0,t)) \times (\Phi^n(v_2,t) - \Phi^n(v_0,t))\|\Big)
\end{multline}
and thus
\begin{multline}
\big| (|K(t)|- |\bar K(t)| )_t\big| \\ \leq  \frac{1}{hy}
 \Big(\int_{\bar K(t_n)} \big\|\partial_t \partial_{v_1} \Phi^n_c(x,t) \times \partial_{v_2} \Phi^n_c(x,t)  
+  \partial_{v_1} \Phi^n_c(x,t) \times \partial_t\partial_{v_2} \Phi^n_c(x,t) \\
- \partial_t(\Phi^n(v_1,t) - \Phi^n(v_0,t)) \times (\Phi^n(v_2,t) - \Phi^n(v_0,t))\\
- (\Phi^n(v_1,t) - \Phi^n(v_0,t)) \times \partial_t(\Phi^n(v_2,t) - \Phi^n(v_0,t))
  \big\|\, d\bar K(t_n)\Big).
\end{multline}
Using the mean value theorem this implies
\begin{multline}\label{Kqd}
\Big|2 (|K(t)|- |\bar K(t)| )_t\Big| \\ \leq  \big\|\partial_t \partial_{v_1} \Phi^n_c(\xi_1,t) \times \partial_{v_2} \Phi^n_c(\xi_1,t)  
+  \partial_{v_1} \Phi^n_c(\xi_2,t) \times \partial_t\partial_{v_2} \Phi^n_c(\xi_2,t)  \\
- \partial_{v_1} \partial_t\Phi^n(\xi_3,t) \times \partial_{v_2} \Phi^n(\xi_4,t) 
- \partial_{v_1} \Phi^n(\xi_5,t) \times\partial_{v_2} \partial_t\Phi^n(\xi_6,t)
  \big\|
\end{multline}
for some $\xi_1,\dots \xi_6 \in \bar K(t_n).$
We have
\begin{equation}\label{xi}
( \partial_t \partial_{x_1} \Phi_c^n)(\xi_1,t) = D( \partial_t \Phi^n)( c(\xi_1),t) \partial_{x_1} c(\xi_1) 
= (\partial_{x_1} \partial_t \Phi^n)( \xi_3,t) + \mathcal{O}(h)
\end{equation}
because of \eqref{eq:partialc} and the regularity of $\Phi^n$.
This, and a straightforward estimate for the second factor in the vector product, leads to
\begin{equation} 
 \partial_t \partial_{v_1} \Phi^n_c(\xi_1,t) \times \partial_{v_2} \Phi^n_c(\xi_1,t)  - \partial_{v_1}\partial_t\Phi^n(\xi_3,t) \times \partial_{v_2} \Phi^n(\xi_4,t) 
= \mathcal{O}(h^3).
\end{equation}
 Using a similar estimate for the remaining terms in \eqref{Kqd} we find
\begin{equation}
 (|K(t)|- |\bar K(t)| )_t = \mathcal{O}(h^3)
\end{equation}
which implies \eqref{eq:app1}$_1$ because of \eqref{ineq:angle-assumption}.\qed

\bibliographystyle{plain}      
\bibliography{bib}

\def\cprime{$'$}
\begin{thebibliography}{10}

\bibitem{AB09}
A.~Alke and D.~Bothe.
\newblock 3d numerical modeling of soluble surfactant at fluidic interfaces
  based on the volume-of-fluid method.
\newblock {\em Fluid Dynamics \& Materials Processing}, 5(4):345--372, 2009.

\bibitem{ABL05}
P.~Amorim, M.~Ben-Artzi, and P.~G. LeFloch.
\newblock Hyperbolic conservation laws on manifolds: total variation estimates
  and the finite volume method.
\newblock {\em Methods Appl. Anal.}, 12(3):291--323, 2005.

\bibitem{ALN11}
P.~Amorim, P.~G. LeFloch, and W.~Neves.
\newblock A geometric approach to error estimates for conservation laws posed
  on a spacetime.
\newblock {\em Nonlinear Anal.}, 74(15):4898--4917, 2011.

\bibitem{BFL09}
M.~Ben-Artzi, J.~Falcovitz, and P.~G. LeFloch.
\newblock Hyperbolic conservation laws on the sphere. {A} geometry-compatible
  finite volume scheme.
\newblock {\em J. Comput. Phys.}, 228(16):5650--5668, 2009.

\bibitem{BL06}
M.~Ben-Artzi and P.~G. LeFloch.
\newblock Well-posedness theory for geometry-compatible hyperbolic conservation
  laws on manifolds.
\newblock {\em Ann. Inst. H. Poincar\'e Anal. Non Lin\'eaire}, 24(6):989--1008,
  2007.

\bibitem{BS10}
M.~R. {Booty} and M.~{Siegel}.
\newblock {A hybrid numerical method for interfacial fluid flow with soluble
  surfactant}.
\newblock {\em J. Comput. Phys.}, 229:3864--3883, 2010.

\bibitem{BPS05}
D.~Bothe, J.~Pr{\"u}ss, and G.~Simonett.
\newblock Well-posedness of a two-phase flow with soluble surfactant.
\newblock In {\em Nonlinear elliptic and parabolic problems}, volume~64 of {\em
  Progr. Nonlinear Differential Equations Appl.}, pages 37--61. Birkh\"auser,
  Basel, 2005.

\bibitem{CHL06}
D.~A. Calhoun, C.~Helzel, and R.~J. LeVeque.
\newblock Logically rectangular grids and finite volume methods for {PDE}s in
  circular and spherical domains.
\newblock {\em SIAM Rev.}, 50(4):723--752, 2008.

\bibitem{CCL94}
B.~Cockburn, F.~Coquel, and P.~G. LeFloch.
\newblock An error estimate for finite volume methods for multidimensional
  conservation laws.
\newblock {\em Math. Comp.}, 63(207):77--103, 1994.

\bibitem{KN}
A.~Dedner, R.~Kl\"ofkorn, and M.~Nolte.
\newblock {DUNE}-{A}lu{G}rid -- {A} parallel-adaptive unstructured grid
  implementation for {DUNE}, in preparation.

\bibitem{DKNO10}
A.~Dedner, R.~Kl{\"o}fkorn, M.~Nolte, and M.~Ohlberger.
\newblock A generic interface for parallel and adaptive discretization schemes:
  abstraction principles and the {DUNE}-{FEM} module.
\newblock {\em Computing}, 90(3-4):165--196, 2010.

\bibitem{Dem09}
A.~Demlow.
\newblock Higher-order finite element methods and pointwise error estimates for
  elliptic problems on surfaces.
\newblock {\em SIAM J. Numer. Anal.}, 47(2):805--827, 2009.

\bibitem{DE07.2}
G.~Dziuk and C.~M. Elliott.
\newblock Finite elements on evolving surfaces.
\newblock {\em IMA J. Numer. Anal.}, 27(2):262--292, 2007.

\bibitem{DE07}
G.~Dziuk and C.~M. Elliott.
\newblock Surface finite elements for parabolic equations.
\newblock {\em J. Comput. Math.}, 25(4):385--407, 2007.

\bibitem{DKM13}
G.~Dziuk, D.~Kr{\"o}ner, and T.~M{\"u}ller.
\newblock Scalar conservation laws on moving hypersurfaces.
\newblock {\em Interfaces Free Bound.}, 15(2):203--236, 2013.

\bibitem{Gie09}
J.~Giesselmann.
\newblock A convergence result for finite volume schemes on {R}iemannian
  manifolds.
\newblock {\em M2AN Math. Model. Numer. Anal.}, 43(5):929--955, 2009.

\bibitem{GW12}
J.~Giesselmann and M.~Wiebe.
\newblock Finite volume schemes for balance laws on time-dependent surfaces.
\newblock In {\em Numerical Methods for Hyperbolic Equations}, pages 251--258.
  CRC Press, London, 2012.

\bibitem{Gil00}
P.~A. Gilman.
\newblock Magnetohydrodynamic ''shallow-water'' equations for the solar
  tachocline.
\newblock {\em The Astrophysical Journal Letters}, 544(1):L79--L82, 2000.

\bibitem{Gir06}
F.~X. Giraldo.
\newblock High-order triangle-based discontinuous galerkin methods for
  hyperbolic equations on a rotating sphere.
\newblock {\em J. Comput. Phys.}, 214(2):447--465, 2006.

\bibitem{JL04}
J.~James and J.~Lowengrub.
\newblock A surfactant-conserving volume-of-fluid method for interfacial flows
  with insoluble surfactant.
\newblock {\em J. Comput. Phys.}, 201(2):685 -- 722, 2004.

\bibitem{LO08.2}
P.~G. LeFloch and B.~Okutmustur.
\newblock Hyperbolic conservation laws on spacetimes. {A} finite volume scheme
  based on differential forms.
\newblock {\em Far East J. Math. Sci. (FJMS)}, 31(1):49--83, 2008.

\bibitem{LON09}
P.~G. LeFloch, B.~Okutmustur, and W.~Neves.
\newblock Hyperbolic conservation laws on manifolds. {A}n error estimate for
  finite volume schemes.
\newblock {\em Acta Math. Sin. (Engl. Ser.)}, 25(7):1041--1066, 2009.

\bibitem{LM13}
D.~Lengeler and T.~Müller.
\newblock Scalar conservation laws on constant and time-dependent riemannian
  manifolds.
\newblock {\em J. Differential Equations}, 254(4):1705 -- 1727, 2013.

\bibitem{LNR11}
M.~Lenz, S.~F. Nemadjieu, and M.~Rumpf.
\newblock A convergent finite volume scheme for diffusion on evolving surfaces.
\newblock {\em SIAM J. Numer. Anal.}, 49(1):15--37, 2011.

\bibitem{RS05}
E.~Reister and U.~Seifert.
\newblock Lateral diffusion of a protein on a fluctuating membrane.
\newblock {\em EPL (Europhysics Letters)}, 71(5):859, 2005.

\bibitem{Ros06}
J.~A. Rossmanith.
\newblock A wave propagation algorithm for hyperbolic systems on the sphere.
\newblock {\em J. Comput. Phys.}, 213(2):629--658, 2006.

\bibitem{SBG01}
D.~A. Schecter, J.~F. Boyd, and P.~A. Gilman.
\newblock ''shallow-water'' magnetohydrodynamic waves in the solar tachocline.
\newblock {\em The Astrophysical Journal Letters}, 551(2):L185--L188, 2001.

\bibitem{Sto90}
H.~A. Stone.
\newblock A simple derivation of the time-dependent convective-diffusion
  equation for surfactant transport along a deforming interface.
\newblock {\em Physics of Fluids A: Fluid Dynamics}, 2(1):111--112, 1990.

\bibitem{OS83}
A.~Van~Oosterom and J.~Strackee.
\newblock The solid angle of a plane triangle.
\newblock {\em Biomedical Engineering, IEEE Transactions on}, BME-30(2):125
  --126, 1983.

\bibitem{WDHJS92}
D.~L. Williamson, J.~B. Drake, J.~J. Hack, R.~Jakob, and P.~N. Swarztrauber.
\newblock A standard test set for numerical approximations to the shallow water
  equations in spherical geometry.
\newblock {\em J. Comput. Phys.}, 102(1):211--224, 1992.

\end{thebibliography}
\end{document}